\documentclass[aap,seceqn,citesort,MSNbibl,dvips]{arximspdf}
\usepackage{subenv}
\usepackage{multicol}

\doi{10.1214/09-AAP633}
\volume{20}
\issue{2}
\pubyear{2010}
\firstpage{660}
\lastpage{695}

\makeatletter

\newcommand{\llangle}{\langle\!\langle}
\newcommand{\rrangle}{\rangle\!\rangle}

\DeclareMathAlphabet\mathcaligr{OMS}{cmsy}{m}{n}

\newtheorem{theorem}{Theorem}
\newtheorem{proposition}{Proposition}
\newtheorem{corollary}{Corollary}
\newtheorem{lemma}{Lemma}

\newcommand{\interleave}{\Vert\hspace*{-1.4pt}\vert}

\makeatother

\begin{document}
\begin{frontmatter}

\title{Brownian coagulation and a version of Smoluchowski's equation on
the circle}
\runtitle{Brownian coagulation on the circle}

\begin{aug}
\author[A]{\fnms{In\'{e}s} \snm{Armend\'{a}riz}\corref{}\ead[label=e1]{iarmendariz@udesa.edu.ar}}
\runauthor{I. Armend\'{a}riz}
\affiliation{Universidad de San Andr\'{e}s}
\address[A]{Departamento de Matem\'{a}tica\\
Vito Dumas 284\\
Victoria\\
B1644BID Buenos Aires\\
Argentina\\
\printead{e1}} 
\end{aug}

\pdfauthor{Ines Armendariz}

\received{\smonth{5} \syear{2006}}
\revised{\smonth{7} \syear{2009}}

%
\begin{abstract}
We introduce a one-dimensional stochastic system where particles
perform independent diffusions and interact through pairwise coagulation
events, which occur at a nontrivial rate upon collision. Under
appropriate conditions
on the diffusion coefficients, the coagulation rates and the initial
distribution of particles, we derive a spatially inhomogeneous
version of the mass flow equation as the particle
number tends to infinity. The mass flow equation is in one-to-one
correspondence with Smoluchowski's coagulation equation. We prove
uniqueness for this equation in a
broad class of solutions, to which the weak limit of the stochastic
system is shown to belong.
\end{abstract}

%

\end{frontmatter}

\section{Introduction}\label{sec1}

Coagulation models describe the dynamics of cluster growth. Particles
carrying
different masses move freely through space, and every time any two of them
get sufficiently close there is some chance that they coagulate into a
single particle, which will be charged with
the sum of the masses of the original pair.

In 1916, Smoluchowski \cite{Smoluchowski}
considered the model
of Brownian particles moving independently in three-dimensional space,
such that any pair coagulates into one particle upon collision. He
derived a system of equations, known as Smoluchowski's coagulation equations,
that describes
the time evolution of the average concentration $\mu_t(m)$
of particles carrying a given mass $m=1, 2, \ldots.$ In this original work,
Smoluchowski ignored the effect of
spatial fluctuations in the mass concentrations,
the equations we write below are thus a natural extension
allowing diffusion in the space variables:
\begin{eqnarray*}
\dot{\mu}_t(m)&=&\frac{1}{2} a(m) \Delta_x\mu_t(x,m)
+
\frac{1}{2}\sum_{m'+m''=m} \kappa(m',m'') \mu_t(x,m') \mu
_t(x,m'')\\
&&{} - \mu_t(x,m) \sum_{m'} \kappa(m,m') \mu_t(x,m').
\end{eqnarray*}
The differentiated term on the right
describes the free motion of a particle with
attached mass $m$ as a Brownian motion with diffusivity rate $a(m)$.
The kernel $\kappa$ is determined by physical considerations, it regulates
the intensity of the coagulation dynamics.
The first sum
corresponds to the increase in the concentration resulting from the coagulation
of two particles whose masses add up to $m$.
The second sum reflects the decrease caused by the coalescence of a particle
carrying mass $m$ with any other particle in
the system. Coagulation phenomena have been studied in many fields of
the applied science; we refer to Aldous's review \cite{Aldousreview}
for a comprehensive survey of the literature.

This paper is concerned with the approximation of Smoluchowski's
equations by stochastic particle models. Concretely, we
are interested in identifying the solution to the coagulation
equations as the mass density
of a system of interacting particles, when the
particle number tends to infinity. This problem has been much studied
in the
spatially homogeneous case, both for discrete and continuous mass
distributions, and with different
choices of coagulation kernel $\kappa$
(cf. \cite{Jeon,NorrisSmol,EibeckWagner,FournierGiet} and references therein).
The relevant stochastic process for these models is the
Marcus and Lushnikov process \cite{Marcus,Lushnikov}; this is
the pure jump
Markov process where
clusters
of size $m$ and $m'$ coagulate into a single cluster of size $m+m'$
at rate $\kappa(m,m')$.

In the spatially inhomogeneous case, on the other hand,
the coagulation mechanism
is highly dependent on the relative position of the particles, hence
the space dynamics plays a predominant role in the particle
interactions of the stochastic system.
In
the original problem proposed by Smoluchowski, for instance,
pairwise collisions and
the ensuing coagulation events are completely determined by the
Brownian paths. The first result for this model
was obtained in 1980 by
Lang and Xanh \cite{LangXanh}
for the case of discrete mass and constant $a, \kappa$, in the limit of
constant mean free time. No progress was made until the forthcoming paper
\cite{NorrisBcoag}, where Norris
proves convergence for
both discrete and continuous mass distributions, and variable
coefficients in a class that
includes the Brownian case.

Over the past few years there has been considerable interest in spatial
models with stochastic dynamics of coagulation. In these models,
particles coagulate at some rate while they remain at less than a
prescribed distance. Deaconu and Fournier \cite{DeaconuFournier}
consider the case when this distance is independent of the particle
number, and let it go to zero after taking the weak limit. The moderate
limit, where the range of interaction is long in the microscopic scale,
is studied by Grosskinsky, Klingenberg and Oelschl\"{a}ger in
\cite{Grosskinsky} in a regime where the dominating particle
interaction are shattering collisions. In the articles
\cite{HammondRez2,HammondRez}, Hammond and Rezakhanlou work in the
constant mean free time limit, for dimensions $d \ge2$.

In this paper we introduce a diffusion model where coagulation occurs
on collision as a result of a random event: $N$
mass-charged
particles perform independent one-dimensional diffusions, and whenever
two particles are at the same location they may coagulate at a
positive rate in their intersection local time. The
new particle is assigned the
sum of the masses of the incoming particles, and the process
continues.

This model is motivated by the problem of establishing the large-scale
dynamics of a system of Brownian particles confined to a thin tube,
interacting through pairwise coagulation when any two of them get
close enough. It would be interesting to determine whether, under
proper scaling of the tube and particle radii in terms of the system
size $N$, the higher-dimensional model can be replaced by the simpler
one-dimensional one.

Note that in one dimension the
problem of instantaneous coagulation on collision
is not interesting: due to the recurrence properties of the
one-dimensional Brownian path, in the limit
we would instantaneously see the distribution of the total mass
in the system
among clusters of macroscopic size. In fact, in order to keep the
average time $T$ during which a tagged particle does not undergo a
collision constant, it is necessary to set coagulation rates which are
inversely proportional to the particle number $N$. The model under
consideration is thus of the constant mean free time type, and in this
sense it is a one-dimensional version of the models studied in
\cite{LangXanh,HammondRez2,HammondRez} and \cite{NorrisBcoag}.

We treat the case where the mass dependent diffusivity rates blow up as
the mass goes to zero, combined with our choice of rates, this leads to
a large-scale model favoring coagulation of large and small particles.

We describe the particle system and state the main results of the paper
in Section~\ref{sec2}. Section \ref{sec3} contains a compactness
result. The next step in the analysis is to prove convergence to a
hydrodynamic limit. This is shown to verify a spatial version of the
mass flow equation, which is closely related to Smoluchowski's
coagulation equation. This is the content of Section \ref{sec4}. In
Section \ref{sec5} we derive a uniqueness result for the solutions to a
broad family of such equations, thereby obtaining a law of large
numbers for the empirical processes of the microscopic model.

\section{Notation and results}\label{sec2}


Consider a positive integer $N$. Let $\mathbb{T}$ stand for the
one-dimensional torus, and $\mathbb{R}_+$ for the half line $[0,
\infty)$. Let
\[
P_0^N=P_0^N( dx^1,dm^1; \ldots
; dx^N,dm^N)
\]
be a sequence of measures on $( \mathbb{T}\times{\mathbb{R}}_+ )^N$
which are symmetric on the pairs $(x^i,m^i)$ and supported on $\sum_i
m^i=1$. Denote by $\mathcaligr{M}_1(\mathbb{T}\times{\mathbb{R}}_+)$
the space of probability measures on $\mathbb{T}\times{\mathbb{R}}_+$
endowed with the weak topology, and by
$\mathcaligr{M}_f(\mathbb{T}\times{\mathbb{R}}_+)$ the space of
positive and finite measures on $\mathbb{T}\times\mathbb{R}_+$.

For a complete separable metric space $\mathcaligr{T}$, we will denote
by $D(\mathbb{R}_+,\mathcaligr{T})$ (or $D(I,\mathcaligr{T})$, $I$ an
interval of the real numbers), the set of right continuous functions
with left limits taking values in $\mathcaligr{T}$, endowed with the
Skorokod topology.

\subsection{The particle model}\label{sec21}
Let
\[
\Phi(m,m')\dvtx{\mathbb{R}}_+ \times{\mathbb{R}}_+ \longrightarrow
{\mathbb{R}}_+
\]
be
a nonnegative, symmetric kernel with the property that it
vanishes when either of the coordinates equals $0$.\vspace*{1pt}

Given a point $\{(x^i,m^i)\}$ in $(\mathbb{T}\times\mathbb{R}_+)^N$,
define
\[
(X^i_{T_0},M^i_{T_0})=(x_i,m_i),\qquad 1\le i\le N,
\]
and set $T_0=0$. Let $k\in\mathbb{N}\cup\{0\}$, and suppose that the process
\[
\xi^N_{\cdot}=\{(X^i_{\cdot}, M^i_{\cdot})\}_{1\le i\le N}\in
C ([0,T_k],\mathbb{T}^N ) \times D ([0,T_k], \mathbb
{R}_+^N )
\]
has already\vspace*{1pt} been defined up to the time $T_k$, a stopping time with
respect to the $\sigma$-algebra $\mathcaligr{F}_t=\sigma\{\xi^N_s, 0\le
s\le t\}$ generated by $\xi^N_{\cdot}$. Consider then a family of $N$
independent Brownian motions $\{\beta^{k,i}_{ \cdot}\} $ on
$\mathbb{T}$ with corresponding diffusion coefficients $a(NM^i_{T_k})$
and initial positions $\beta^{k,i}_0=X^i_{T_k}$. For each pair $i<j$,
denote by $L^{k,ij}$ the intersection local time of the $i$th and $j$th
particles; that is, the local time at the origin of the difference
$\beta^{k,i}-\beta^{k,j}$. Let $\{\epsilon^{k,ij}\}_{1\le i<j\le N}$ be
a sequence of ${N\choose2}$ independent, parameter one-exponential
random variables, and define the stopping times
\[
\tilde{T}_{k+1}=\min_{i<j} \{ T^{ij}_{k+1} \},\qquad
T^{ij}_{k+1}=\inf\biggl\{t\ge0,
\frac{\Phi(NM^i_{T_k},NM^j_{T_k})}{N} L^{k,ij}_t > \epsilon
^{k,ij} \biggr\}.
\]
Let then $T_{k+1}=T_k+\tilde{T}_{k+1}$, and for $1\le i\le N$, set
\begin{eqnarray*}
X^i_t &=& \beta^{k,i}_{t-T_k},\qquad T_k < t \le T_{k+1}, \\
M^i_t &=& M^i_{T_k},\qquad T_k< t < T_{k+1},
\\
M^i_{T_{k+1}} &=& \cases{
M^i_{T_k}+M^j_{T_k}, &\quad if
$\tilde{T}_{k+1}=T^{ij}_{k+1}$, for some $j>i$;\cr
0, &\quad if
$\tilde{T}_{k+1}=T^{ji}_{k+1}$, for some
$j<i$;\cr
M^i_{T_k}, &\quad otherwise.}
\end{eqnarray*}

We will denote by $L^{ij}$ the intersection local time of the $i$th
and $j$th particles $X^i$ and $X^j$. Note that between two consecutive stopping
times $T_k$ and $T_{k+1}$ the identity $L^{ij}=L^{k,ij}$ holds.

The dynamics is well defined except for those configurations where two
or more coagulation events occur simultaneously. Let us briefly show
that the set of such configurations has measure zero and can therefore
be neglected.\vspace*{1pt} We first show that this is the case on the time interval
$[0,T_1]$, when $X^i=\beta^{1,i}, 1\le i\le N$. The $i$th and $j$th
masses will coagulate at a time belonging to the support of the
measure~$dL^{ij}$, which equals the zero set of
$\beta^{1,i}-\beta^{1,j}$. Fix two pairs of indices $i<j$ and $k<l$. If
the four indices are different, then $\beta^{1,i}-\beta^{1,j}$ and
$\beta ^{1,k}-\beta^{1,l}$ perform independent diffusions, and from the
fact that point sets are polar for Brownian motion in $2$ or higher
dimensions it follows that these diffusions do no vanish at the same
time. Let there be a repeated index: $j=k$, say. Set $\alpha=\frac
{a(NM^j)+a(NM^l)}{a(NM^j)}$. Then
\[
U_t=\alpha[\beta^{1,i}_t-\beta^{1,j}_t]+\beta^{1,j}_t-\beta^{1,l}_t
\quad\mbox{and}\quad V_t=\beta^{1,j}_t-\beta^{1,l}_t
\]
are independent diffusions in $\mathbb{T}^2$. The previous argument
implies that
with probability 1 they never vanish simultaneously, and hence the
same applies to $\beta^{1,i}-\beta^{1,j}$ and $\beta^{1,j}-\beta
^{1,l}$. We conclude
that in any case
the zero sets of
$\beta^{1,i}-\beta^{1,j}$ and $\beta^{1,k}-\beta^{1,l}$ are
disjoint, with
probability $1$. By repeating the argument on each interval
$[T_k,T_{k+1}], k\ge1$, it follows that outside a set of measure zero
there are no conflicting coagulation
events.

\subsection{Martingales}\label{sec22}
Let $P^N$ be the measure on $C(\mathbb{R}_+, {\mathbb{T}}^N) \times
D(\mathbb{R}_+,{\mathbb {R}}_+^N)$ determined by the process
$\xi_{\cdot}^N$. There is the representation
\[
M^i_t=m^i+\int_0^t \sum_{i<j\le N} M^j \,dE^{ij}-
\int_0^t \sum_{1\le k<i} M^i \,dE^{ki},
\]
where $dE^{ij}$ is a counting measure, $E^{ij}([0,t])=0,1$ with
$P^N$-probability $1$, depending on whether the $i$th and $j$th
particles have coagulated by time $t>0$, and
\[
\mathcaligr{H}^{ij}_t=E^{ij}([0,t])-\int_0^t \frac{\Phi
(NM^i_s,NM^j_s)}{N}\,
dL^{ij}
\]
is a martingale.

In general, if $f$ is a bounded function on $(\mathbb{T}\times\mathbb
{R}_+)^N$ with two continuous, bounded derivatives in the space
coordinates, then
\begin{eqnarray*}
&&f(\xi_t)-\int_0^t \sum_{i=1}^N a(NM^i_s) \,\frac{\partial^2
f}{\partial{X^i}^2}(\xi_s) \,ds \\
&&\qquad{} -
\int_0^t \sum_{i<j} [f(\xi^{ij}_s)-f(\xi_s) ] \frac
{\Phi
(NM_s^i, NM_s^j)}{N} \,dL^{ij}
\end{eqnarray*}
is an $(\mathcaligr{F}_t, P)$-martingale. Given $\xi=\{(X^i,M^i)\}$,
$\xi^{ij}$ here is given by
\[
(\xi^{ij})^k= \cases{
(X^i, M^i+M^j), &\quad if $k=i$;\cr
(X^j, 0), &\quad if $k=j$;\cr
(X^k, M^k), &\quad otherwise.}
\]

\subsection{Scaling}\label{sec23}

The reason for the choice of scaling of the coagulation rate is quite
straighforward: if a hydrodynamic description is to hold, then it is
necessary that $O(N)$ mass charged particles remain in the system at
all times (note that although total mass is conserved, the number of
particles carrying positive mass decreases by one after each
coagulation event). Therefore a generic particle will see some fraction
of the other $O(N)$ mass charged particles over any fixed time
interval, while it would still be expected to coagulate with only
$O(1)$ of them. This forces the rate $\Phi$ to be typically of order
$1/N$.

\subsection{Assumptions}\label{sec24}

We will consider the mapping
\[
\Pi_N\dvtx C (\mathbb{R}_+, {\mathbb{T}}^N ) \times D
(\mathbb{R}_+,\mathbb{R}_+^N
)\rightarrow
D \bigl(\mathbb{R}_+,\mathcaligr{M}_1 (\mathbb{T}
\times\mathbb{R}_+) \bigr)
\]
such that
\[
\Pi_N (\{X^i_\cdot,M^i_\cdot\} )=\sum_{1 \le i \le N}
M^i_\cdot
\delta_{(X^i_\cdot, NM^i_\cdot)}
\]
and denote by $Q^N$ the measure on $D
(\mathbb{R}_+,\mathcaligr{M}_1(\mathbb{T}\times\mathbb {R}_+) )$
induced by $\Pi_N$,
\[
Q^N=P^N \circ\Pi_N^{-1}.
\]

In order to derive a hydrodynamic limit for $Q^N$, we
need to specify some technical conditions on the
coalescing kernel $\Phi$, the diffusion coefficients $a$,
the initial measure
$P_0^N$ and the profile $\nu$.

The kernel $\Phi$ satisfies
a Lipschitz condition away from the origin: for each $L>0$
there exists a positive constant $\Gamma(L)$ such that
%
%
\begin{subequation}
\begin{equation}\label{equ21a}
|\Phi(m+m'',m')-\Phi(m,m') | \le
\Gamma(L) m'' \qquad\mbox{whenever } m>L .
\end{equation}
It will also be assumed that
there exists $0 \le\mathfrak{p}\le1/2$ such that
\begin{equation}\label{equ21b}
\Phi(m,m')\le c (m^{\mathfrak{p}}+{m'}^{\mathfrak{p}})
1_{[m>0, m'>0]} 
\end{equation}
\end{subequation}
for some positive constant $c$.

We assume that the mapping
$a\dvtx(0,\infty)\rightarrow(0,\infty)$ is nonincreasing, so that
particles diffuse at a slower rate as they gain mass. As a consequence
the kernel $\Phi$ and the diffusion coefficients do not grow
simultaneously. We set $a(0)=0$.

We are ready to introduce the coagulation rates $\kappa\dvtx\mathbb{R}_+
\to\mathbb{R}
_+$ appearing in the hydrodynamic equation,
\[
\kappa(m,m')=\Phi(m,m') [a(m)+a(m') ].
\]
In
order to study convergence and derive the uniqueness of the
limit it will be useful to consider
\[
\omega(m)=[1+c+a(1)] [m^{\mathfrak{p}}+a(m)+1 ],
\]
it verifies $\kappa(m,m')\le\omega(m) \omega(m')$. We will then
require that
%
%
\begin{equation}\label{equ22}
a(m)^{-1/2} \omega(m) \mbox{ be a subadditive function of $m$.} %
\end{equation}

Conditions (\ref{equ21a}), (\ref{equ21b}) and (\ref{equ22}) are for
instance satisfied by
\[
\Phi(m,\tilde{m})=C(m^\alpha+\tilde{m}^\alpha) \quad\mbox
{and}\quad a(m)=\frac{1}{m^\beta}1_{[m>0]} \qquad\mbox{with }\alpha\le
\frac{1}{2} \mbox{ and }\beta\le1.
\]
In this case
$\kappa(m, \tilde{m})=C(m^\alpha+\tilde{m}^\alpha) (\frac
{1}{m^\beta}+\frac{1}{\tilde{m}^\beta} )$.\vspace*{1pt}

We will assume that there exists an initial profile
$\nu\in\mathcaligr{M}_1({\mathbb {T}\times {\mathbb{R}}_+})$ such that
the empirical distributions $\sum_i m^i\delta_{(x^i,Nm^i)}$ converge in
distribution to $\delta_{\nu}$ as $N\rightarrow\infty$, where by
$\delta_{(x^i,Nm^i)}$ (resp., $\delta_{\nu}$) we denote the Dirac
measure with a unit atom at $(x^i,Nm^i)$ (resp., $\nu$).

The initial measures $P_0^N$ will satisfy
%
%
\begin{equation}\label{equ23}
E^{P_0^N} \Bigl[N\sum(m^i)^2 \Bigr]<C \quad\mbox{and}\quad
E^{P_0^N} \biggl[\frac{1}{N} \sum a(N m^i)^2 \biggr]<C'
\end{equation}
for some constants $C,C'>0$, uniformly in $N$. In particular both
$\langle m,\nu\rangle$ and $\langle a(m)^2/m, \nu\rangle$ are finite.
In fact, the following assumption will hold: there exists a finite
measure $\nu^*(dm)$ such that
%
%
\begin{equation}\label{equ24}
\nu(dx,dm) \le\nu^*(dm) \,dx \qquad\mbox{with }
\biggl\langle m + \frac{a(m)^2}{m}, \nu^* \biggr\rangle< \infty.
\end{equation}

\subsection{Results}\label{sec25}

The first theorem of the paper is a tightness result:
\begin{theorem}\label{theo1}
Assume that (\ref{equ21a}), (\ref{equ21b}), (\ref{equ22}) and (\ref
{equ23}) hold. Then the sequence of measures $Q^N$ on $D
(\mathbb{R}_+,\mathcaligr{M}_1(\mathbb{T}\times\mathbb {R}_+) )$ is
relatively compact, and all limit points are concentrated on continuous
paths.
\end{theorem}

The next two results concern the properties satisfied by any weak limit of
the empirical distributions as we pass to the limit in the particle number.
The first result provides some estimates that will ensure
the hydrodynamic equation is well defined, then Theorem \ref{theo2}
identifies this equation, thereby establishing an existence
result.

Given a kernel admitting a representation $\mu(x,dm) \,dx \in
\mathcaligr{M}_f(\mathbb{T}\times\mathbb{R}_+)$ and a bounded test
function $f(m)$, we will denote by $\llangle f,\mu\rrangle$ the single integral
$\int_{\mathbb{R}_+} f(m) \mu(x,dm)$. This clearly determines a signed
measure on $\mathbb{T}$ by
\[
\int_{\mathbb{T}} h(x)\llangle f,\mu\rrangle \,dx=\int_{\mathbb{T}\times
\mathbb{R}_+}
h(x) f(m) \mu(x,dm) \,dx.
\]
\begin{proposition}\label{prop1}
Let $Q$ be a weak limit of the sequence $Q^N$ of probability measures
on $C (\mathbb{R}_+,\mathcaligr{M}_1(\mathbb{T}\times\mathbb {R}_+) )$.
Then $Q$ is supported on the set of paths $\mu_{\cdot}(dx,dm)$ whose
marginal $\mu_t(dx,\mathbb{R}_+) \ll dx$ on $\mathbb{T}$ for all $t$,
$\mu_t(dx,dm)=\upsilon_t(x,dm) \,dx$. Moreover, the following
inequalities hold with $Q$-probability 1:
%
%
\begin{equation}\label{equ25}
\sup_{t\ge0}
\biggl\| \biggl\langle\!\biggl\langle\frac{\omega(m)}{m},
\upsilon_t \biggr\rangle\!\biggr\rangle\biggr\|_{\infty}
<\infty
\end{equation}
and
%
%
\begin{equation}\label{equ26}
\sup_{t\le T} \langle m,\mu_t\rangle<\infty
\end{equation}
for any fixed final time $T$.
\end{proposition}

Denote by $C_b^2(\mathbb{T}\times\mathbb{R}_+)$ the space of
continuous, bounded functions
on $\mathbb{T}\times\mathbb{R}_+$ which have continuous, bounded
derivatives in the space
variable up to the second order.
\begin{theorem}\label{theo2}
Let $Q$ be a weak limit of the sequence $Q^N$, as in Proposition \ref
{prop1}, and consider
$f$ in $C^{2}_b(\mathbb{T}\times\mathbb{R}_+)$. Then, with
$Q$-probability 1, a path
$\mu_s(dx, dm)$ satisfies
%
%
\begin{eqnarray}\label{equ27}\hspace*{30pt}
\langle f,\mu_t \rangle-\langle f,\nu\rangle
&=&
\int_0^t \biggl\langle\frac{1}{2} a(m) \,\frac{\partial^2
f}{\partial x^2},
\mu_s \biggr\rangle \,ds \nonumber\\
&&{} + \int_0^t \int_{\mathbb{T}} \int_{\mathbb
{R}_+\times\mathbb{R}_+}
\frac{ [f(x,m+m')-f(x,m) ]}{m'} \kappa(m,m') \\
&&\hspace*{75pt}{}
\times\upsilon_s(x,dm) \upsilon_s(x,dm') \,dx \,ds,\nonumber
\end{eqnarray}
if $t\ge0$. In this equation $\nu$ is the initial profile of the model,
and the kernel $\upsilon_s(x,dm)$ is such that
\[
\mu_s(dx,dm)=\upsilon_s(x,dm) \,dx \qquad\mbox{for
all $s\ge0$, $Q$-a.e.}
\]
\end{theorem}

Equation (\ref{equ27}) describes the evolution in time of the mass
flow: if we
decompose its solutions
as $\upsilon_s(x,dm) \,dx=m \hat{\upsilon}_s(x,dm) \,dx$, an elementary
computation proves that $\hat{\upsilon} \,dx$ satisfies Smoluchowski's
coagulation equation with kernel $\kappa$. Theorem \ref{theo1} asserts
that all weak
limits of the measures $Q^N$ are supported on $\mathcaligr
{M}_1(\mathbb{T}
\times\mathbb{R}_+)$;
in terms of the concentration densities $\hat{\upsilon}_s(x,dm)$, this
means that mass is conserved, or equivalently, that there is no gelation
phenomenon.

The method
applied to derive
these results relies heavily on stochastic calculus
computations, we try to make these quite detailed in the proof of
Theorem \ref{theo1}
and
give only an outline later on.

In \cite{NorrisBcoag}, Norris introduces a
method for proving existence and uniqueness for a general
class of $d$-dimensional diffusion models with coagulation; briefly
put, this consists on
approximating the corresponding version of (\ref{equ27}) in his paper
by a
system that depends on the
coalescing kernel $\kappa$ only through its values on a given compact
set. In Section \ref{sec5} we develop
a simplified version of his technique to obtain a
uniqueness result for the solutions of a broad family of mass flow equations.
Section \ref{sec5} may be read independently of the rest of the paper.

Some brief consideration shows that the right-hand side of
(\ref{equ27}) is well defined
in a proper subset
of $C (\mathbb{R}_+,\mathcaligr{M}_1(\mathbb{T}\times\mathbb
{R}_+) )$ consisting
of those paths $\eta$ whose
marginal $\eta_t(dx, \mathbb{R}_+)$ has a
density with respect to Lebesgue measure satisfying some integrability
conditions.
It is easy to see that in fact (\ref{equ25}) and (\ref{equ26}) are
enough, and then
Proposition \ref{prop1} says that
all weak limits of the sequence $Q^N$ are supported on configurations
where the right-hand side
of (\ref{equ27}) can be evaluated. This observation motivates the following
definition: we will denote by
$\mathcaligr{D}(\omega)$ the subset of $C(\mathbb{R}_+,\mathcaligr{
M}_f(\mathbb{T}\times\mathbb{R}_+))$ of those paths $\eta$ whose marginal
$\eta_t(dx,\mathbb{R}_+) \ll dx$, and such
that
(\ref{equ25}) and (\ref{equ26}) are satisfied. Note that the map $\omega
$ depends on the
diffusivity $a$ and the coagulation rate $\Phi$, and then so does
$\mathcaligr{D}(\omega)$.

As a particular case of Theorem \ref{theo3} in Section \ref{sec5}, we have:
\begin{corollary}\label{coro1}
Assume conditions (\ref{equ21b}), (\ref{equ22}) and (\ref{equ24}) on
the coagulation rate
$\kappa$
and the initial measure $\nu$, respectively.
Then for any $T\ge0$, (\ref{equ27}) has at most one solution $\{\mu_t\}
_{0\le
t\le T}$
in $\mathcaligr{D}(\omega)$.
\end{corollary}

The four preceding results imply that the sequence of probability
measures $Q^N$ converges to the Dirac measure
concentrated on the unique solution in $\mathcaligr{D}(\omega)$
to~(\ref{equ27}).

\section{Existence of a weak limit}\label{sec3}

In order to simplify notation, we will often omit the dependence of
the masses and positions on the time parameters whenever we think
that this would not
lead to confusion.
For instance, in an integral where time is
parametrized by $s$, $M^i$ and $X^i$ should be read as $M^i_s$
and $X^i_s$, respectively.

Throughout the article, $\Gamma$ will denote a positive constant. Unless
we are particularly interested in keeping track of its growth or
dependence on the parameters, we will use the same letter $\Gamma$ to denote
constants on consecutive lines which may be different, or constants
appearing in totally unrelated computations.

Let us consider a fixed final time $T>0$ for the rest of the paper.
We will prove a version of Theorems \ref{theo1}, \ref{theo2} and
Proposition \ref{prop1} on
the compact interval $[0,T]$; the fact that the value of $T$ is
arbitrary will then
imply that these results hold as stated in the previous section.

The following estimates will be necessary to derive
Theorem \ref{theo1}; we postpone their proofs until the end of this section.
\begin{lemma}\label{lemma1}
There exist nonnegative constants $C(T), C'(T)$ which depend on
the diffusivity $a$, the kernel $\phi$ and the bounds appearing
in (\ref{equ23}),
such that
%
%
\begin{eqnarray}\label{equ31}
E^{P^N} \biggl[ N\sum_i [M^i_T]^2 \biggr] &<& C(T) \quad\mbox{and}\\
\label{equ32}
E^{P^N} \biggl[ \int_0^T \sum_{i<j} M^i M^j \Phi(NM^i,NM^j)
\,dL^{ij} \biggr] &<& C'(T)
\end{eqnarray}
hold uniformly in $N$.
\end{lemma}
\begin{lemma}\label{lemma2}
Given $\epsilon>0$, there is $\delta>0$ such that
\[
\lim_{N\rightarrow\infty}
P^N \biggl[ \mathop{\mathop{\sup}_{0\le s\le t\le T}}_{t-s \le\delta}
\int_s^t \sum_{i<j} M^i M^j \Phi(NM^i,NM^j) \,dL^{ij}> \epsilon
\biggr] < \epsilon.
\]
\end{lemma}

Denote by $C_b(\mathbb{T}\times\mathbb{R}_+)$ the space of bounded,
continuous functions
on $\mathbb{T}\times\mathbb{R}_+$ with the topology determined by
uniform convergence
over compact sets. Let $\{f_k, k \in\mathbb{N}\}$ be a dense,
countable family in $C_b(\mathbb{T}\times\mathbb{R}_+)$. Then the distance
\[
\varrho(\mu,\nu)=\sum_{k \in\mathbb{N}}\frac{1}{2^k}
\frac{|\langle f_k, \mu\rangle-\langle f_k,\nu
\rangle|}{1+|\langle f_k,\mu\rangle- \langle f_k,\nu\rangle|}
\]
defines a metric on $\mathcaligr{M}_1(\mathbb{T}\times\mathbb
{R}_+)$ which is
compatible with the
weak topology.
There is the associated modulus of continuity
\[
\omega_{\mu}(\gamma)=\mathop{\mathop{\sup}_{0\le s\le t\le
T}}_{t-s \le\gamma}
\varrho(\mu_t,\mu_s).
\]
\begin{pf*}{Proof of Theorem \protect\ref{theo1}}
We refer to Chapter 4 in \cite{KipnisLandim} for a presentation of
the Skorokod's
topology as well as the characterization of the relatively compact
sets in $D([0,T],\mathcaligr{M}_1(\mathbb{T}\times\mathbb{R}_+))$.
Note that condition (ii) below implies that, provided
the sequence $Q^N$ has limit points, these will be supported on
$C([0,T],\mathcaligr{M}_1(\mathbb{T}\times\mathbb{R}_+))$.

By a version of Prokhorov's theorem applied to this setting (cf.
\cite{BillingsleyWC}, Chapter 3), the theorem will follow if we can
show that:
\[
\mbox{\quad\hspace*{3pt}(i)}\hspace*{26.1pt} \mbox{ For every $\varepsilon>0$}\qquad
\lim_{M\uparrow\infty} \limsup_{N \rightarrow\infty}
Q^N \biggl[\sup_{0\le t\le T} \mu_t (m>M )>\varepsilon
\biggr]=0\hspace*{14.1pt}
\]
and
\[
\mbox{\quad(ii)}\hspace*{54.5pt} \mbox{ For every $\varepsilon>0$}\qquad\lim_{\gamma
\downarrow0}
\limsup_{N\rightarrow\infty}
Q^N [ \omega_{\mu}(\gamma)>\varepsilon ]=0.\hspace*{42.5pt}
\]
Note that $\langle m,\mu_t \rangle$
is nondecreasing for $0\le t\le T$, $Q^N$-a.e.,
a fact we will repeatedly use in the course of the article.
Then (i) is an easy consequence of
(\ref{equ31}) in Lemma \ref{lemma1} and Chebyshev's inequality.

In order to conclude (ii) it will be enough to prove that given
$f\in C_b^2(\mathbb{T}\times\mathbb{R}_+)$, $f$ Lipschitz in $m$, we
can control
\[
Q^N \biggl[ {\mathop{\mathop{\sup}_{0\le s\le t\le T}}_{t-s \le\gamma}
}|\langle f, \mu_t\rangle- \langle f,\mu_s \rangle|>
\varepsilon\biggr] < \varepsilon
\]
provided $N$ and $\gamma$ are taken to be sufficiently large and small,
respectively.

Define the stopping time
\[
\tau=\inf\biggl\{t\ge0, \max_i(M^i_t) \biggl[1+\sum_j m^j
a(Nm^j) \biggr]>N^{-1/4} \biggr\};
\]
by Chebyshev's inequality we compute
\[
P^N [\tau\le T ]<\frac{[C(T)(1+C')]^{1/2}}{N^{1/4}},
\]
where $C(T)$ and $C'$ are the constants appearing on the right
of (\ref{equ31}) and the second inequality in (\ref{equ23}),
respectively. By stopping the process as soon as $\tau$
is achieved, we may assume that
%
%
\begin{equation}\label{equ33}
P^N \biggl[ \biggl\{\max_i(M^i_t) \biggl[1+\sum_j m^j a(Nm^j) \biggr]\le
N^{-1/4} \biggr\} \biggr]=1.
\end{equation}

Applying It\^{o}'s formula to $f$, we have
\begin{eqnarray*}
&&\langle f,\mu_t \rangle- \langle f, \mu_s \rangle\\
&&\qquad=\int_s^t \sum_{i} M^i
\,\frac{\partial f}{\partial
x}\, (X^i,NM^i) \,dX^i\\
&&\qquad\quad{} +
\frac{1}{2}\int_s^t \sum_i M^i \,\frac{\partial^2 f}{\partial
x^2}\,(X^i,NM^i) a(NM^i) \,ds\\
&&\qquad\quad{} + \int_s^t \sum_i
[F_N(X^i,M^i+M^j)-F_N(X^i,M^i)-F_N(X^i,M^j) ]\,
dE^{ij},
\end{eqnarray*}
where we denote $F_N(x,m)=mf(x,Nm)$.

Let $\gamma>0$. Doob's inequality, (\ref{equ33}) and the monotonicity
of $a$ imply
\begin{eqnarray*}
&& P^N \biggl[\mathop{\mathop{\sup}_{0\le s\le t\le T}}_{t-s \le\gamma}
\biggl|\int_s^t \sum_{i} M^i \,\frac{\partial f}{\partial
x}\, (X^i,NM^i) \,dX^i\biggr|
> \frac{\varepsilon}{3} \biggr]\\
&&\qquad\le P^N \biggl[\sup_{0\le t\le T}
\biggl|\int_0^t \sum_{i} M^i \,\frac{\partial f}{\partial
x}\, (X^i,NM^i) \,dX^i\biggr|
> \frac{\varepsilon}{6}
\biggr]\\
&&\qquad\le\frac{\Gamma(f,\varepsilon,T)}{N^{1/4}},
\end{eqnarray*}
where $\Gamma$ is a positive constant that does not depend on $N$.
By taking $\gamma$ such that
%
%
\begin{equation}\label{equ34}
C' \biggl\|\frac{\partial^2f}{\partial x^2} \biggr\|_{\infty}
\gamma\le\frac{\varepsilon^2}{3}
\end{equation}
we also obtain
\[
P^N \biggl[\mathop{\mathop{\sup}_{0\le s\le t\le T}}_{t-s \le\gamma}
\biggl|\int_s^t \sum_{i} M^i \,\frac{\partial^2 f}{\partial
x^2}\, (X^i,NM^i) a(NM^i) \,ds\biggr|
> \frac{\varepsilon}{3} \biggr]\le\frac{\varepsilon}{3}.
\]

It remains to estimate the Poisson integral
\begin{eqnarray*}
&&\int_s^t \sum_i
[F_N(X^i,M^i+M^j)-F_N(X^i,M^i)-F_N(X^i,M^j) ]
\,dE^{ij}\\
&&\qquad=\mathcaligr{H}_F(0,t)-\mathcaligr{H}_F(0,s)+I_F(s,t),
\end{eqnarray*}
if $I_F(s,t)$ denotes the integral
\[
\int_s^t\sum_{i<j}
[F_N(X^i,M^i+M^j)-F_N(X^i,M^i)-F_N(X^i,M^j) ]
\frac{\Phi(NM^i,NM^j)}{N}\,dL^{ij}
\]
and $\mathcaligr{H}_F(0,t)$ is a martingale collecting the remaining
terms. Its quadratic variation is given by
\[
\int_0^t \sum_{i<j}
[F_N(X^i,M^i+M^j)-F_N(X^i,M^i)-F_N(X^i,M^j) ]^2
\frac{\Phi(NM^i,NM^j)}{N}\,dL^{ij}.
\]
Note that
\[
|F_N(x,m+m')-F_N(x,m)-F_N(x,m')|\le\Gamma(f) [(m+m')\wedge(
Nmm') ].
\]
In particular, due to assumption (\ref{equ33}) on the mass sizes,
(\ref{equ32}) and Doob's inequality, we obtain
\[
P^N \biggl[{\mathop{\mathop{\sup}_{0\le s\le t\le T}}_{t-s \le\gamma}}
|\mathcaligr{H}_F(s,t)|>\frac{\varepsilon}{6} \biggr]\le
\frac{\Gamma(f,\varepsilon,T)}{N^{1/4}},
\]
which will decay to $0$ as we pass to the limit $N\rightarrow\infty$.
Finally,
\[
|I_F(s,t)|\le\Gamma(f)\int_s^t\sum_{i<j}
M^iM^j \Phi(NM^i,NM^j) \,dL^{ij}.
\]
The result now follows by taking
\[
\gamma\le\gamma_1 \wedge\delta,
\]
where $\gamma_1$ satisfies (\ref{equ34}) and $\delta$ is the value
given by Lemma \ref{lemma2} when $\epsilon$ is set equal to
$\varepsilon /6[1+\Gamma(f)]$.
\end{pf*}
\begin{pf*}{Proof of Lemma \protect\ref{lemma1}}
The proof of this lemma will follow from repeated applications of
It\^{o}--Tanaka's theorem; see \cite{IkedaWatanabeSDE} for an
exposition of this and related formulas. We write
%
%
\begin{equation}\label{equ35}\hspace*{25pt}
N\sum_i [M^i_t]^2 = N\sum_i [m^i]^2 + \mathcaligr{H}_t +2\int_0^t
\sum_{i<j} M^i_s M^j_s \Phi(NM^i_s,NM^j_s) \,dL^{ij},
\end{equation}
where $\mathcaligr{H}_t$ is the $P^N$-martingale
\[
\mathcaligr{H}_t=\int_0^t \sum_{i<j}
2NM^i_s M^j_s \biggl[dE^{ij}-\frac{\Phi(NM^i_s,NM^j_s)}{N}\,
dL^{ij} \biggr].
\]
We focus on the last term of (\ref{equ35}). Given $\zeta>0$, let
$g_{\zeta} \in C(\mathbb{T})\cap C^2(\mathbb {T}-\{0\})$ be a positive,
even function that equals $|x|$ in a small interval containing the
origin, vanishes outside $[-1/4,1/4]$ and satisfies
$\sup_{x\in\mathbb{T}} g_{\zeta}(x)\le \zeta$. We then have
%
%
\begin{eqnarray}\label{equ36}
&&\int_0^t\sum_{1\le i < j \le N} M^i M^j \Phi(NM^i,NM^j) \,dL^{ij}
\nonumber\\[-8pt]\\[-8pt]
&&\qquad=A_1(0,t)-A_2(0,t)-A_3(0,t)-A_4(0,t)-A_5(0,t)\nonumber
\end{eqnarray}
with
\begin{eqnarray*}
A_1(0,t) &=& \sum_{1\le i < j \le N} M^i_t M^j_t \Phi(NM^i_t,NM^j_t)
g_{\zeta}(X^i_t-X^j_t)\\
&&{} -\sum_{1\le i< j \le N} m^i m^j \Phi(Nm^i,Nm^j)
g_{\zeta}(X^i-X^j), \\
A_2(0,t) &=& \int_0^t\sum_{1\le i < j \le N} M^i M^j \Phi(NM^i,NM^j)
g_{\zeta}'(X^i-X^j) [dX^i-dX^j ]
\end{eqnarray*}
and
\begin{eqnarray*}
A_3(0,t) &=& \frac{1}{2} \int_0^t \sum_{i<j} M^i M^j\Phi(NM^i,NM^j)
g_{\zeta}''(X^i-X^j)\\
&&\hspace*{37pt}{}\times
[a(NM^i)+a(NM^j) ] \,du .
\end{eqnarray*}
The function $g_{\zeta}''$ appearing in the formula for $A_3(0,t)$ stands
for what is left of the second derivative of $g_{\zeta}$ (in the sense of
distributions) after substracting $2\delta_0$, $\delta_0$ the Dirac
measure at the origin. The terms
$A_4(0,t)$ and $A_5(0,t)$ correspond to the coagulation martingale and
its compensator,
\begin{eqnarray*}
A_4(0,t) &=& \int_0^t \mathop{\mathop{\sum}_{i<k}}_{j}
D_{N}(M^i,M^j,M^k)
g_{\zeta}(X^i-X^j)\\[-6pt]
&&\hspace*{28pt}{} \times\biggl[dE^{ik}-\frac{\Phi(NM^i,NM^k)}{N}
\,dL^{ik} \biggr],\\
A_5(0,t) &=& \int_0^t \mathop{\mathop{\sum}_{i<k}}_{j}
D_{N}(M^i,M^j,M^k)
g_{\zeta}(X^i-X^j) \frac{\Phi(NM^i,NM^k)}{N} \,dL^{ik},
\end{eqnarray*}
where $D_{N}(m,m',m'')$ is defined as
\begin{eqnarray*}
&&(m+m'') m'\Phi(Nm+Nm'',Nm')\\
&&\qquad{} - m m'\Phi(Nm,Nm')-m''m'\Phi(Nm'',Nm') .
\end{eqnarray*}
In deriving the formula for $A_4$, we have used that at the
time when the masses $M^i$ and $M^k$ coagulate,
the $i$th and $k$th particles are occupying the same position.

We will study these terms separately. Replacing
$\Phi(m,m')\le c [m^{\mathfrak{p}}+m'^{\mathfrak{p}} ]$ in
the definition of
$A_1$ gives
%
%
\begin{equation}\label{equ37}
|A_1(0,t)| \le4 c \zeta\sum_{i} M^i_t[N M^i_t]^{\mathfrak{p}}.
\end{equation}
The bounded variation term $A_3(0,t)$ may be similarly controlled,
%
%
\begin{equation}\label{equ38}
|A_3(0,t)|\le\Gamma(\zeta) \int_0^t \biggl[1+\sum_i M^i_s
[NM^i_s]^{\mathfrak{p}} \biggr] \biggl[1+ \sum_i m^i a(Nm^i) \biggr] \,ds .
\end{equation}

In order to bound $A_5$, we first notice that by the Lipschitz
assumption (\ref{equ21a}) on $\Phi$ we have
%
%
\begin{eqnarray}\label{equ39}
&&|\Phi(Nm+Nm'',Nm')-\Phi(Nm,Nm')| \nonumber\\[-8pt]\\[-8pt]
&&\qquad \le
\Gamma(1) Nm'' 1_{\{Nm\ge1\}}
+ ( 1+[Nm']^{\mathfrak{p}}+[Nm'']^{\mathfrak{p}} ) 1_{\{
Nm<1\}} .\nonumber
\end{eqnarray}
It then follows that
%
%
\begin{eqnarray}\label{equ310}\qquad
|A_5(0,t)| & \le & \Gamma\zeta \biggl[\int_0^t \sum_{i<k}
M^iM^k \Phi(NM^i,NM^k) \, dL^{ik} \nonumber\\
&&\hspace*{17.30pt}{} + \biggl(1+\sum_j M^j_t[NM^j_t]^{\mathfrak{p}} \biggr)
\\
&&\hspace*{28.09pt}{}\times \int
_0^t \sum
_{i<k} \frac{1}{N^2} \Phi(NM^i,NM^k) \,dL^{ik} \biggr].\nonumber
\end{eqnarray}
We have
\begin{eqnarray*}
&&\int_0^t \sum_{i<k} \frac{1}{N^2} \Phi(NM^i,NM^k) \,dL^{ik}
\\
&&\qquad=
\sum_i \frac{1}{N} 1_{\{m^i>0 \}} - \sum_i
\frac{1}{N} 1_{\{M^i_t>0 \}}
\\
&&\qquad\quad{}
-\int_0^t \sum_{i<k}
\frac{1}{N} \biggl[dE^{ik}-\frac{\Phi(NM^i,NM^k)}{N} \,dL^{ik} \biggr] .
\end{eqnarray*}
The last term above is a martingale, hence
\[
E^{P^N} \biggl[\int_0^t \sum_{i<k}
\frac{1}{N^2} \Phi(NM^i,NM^k) \,dL^{ik} \biggr]
\le 1
\]
and
\begin{eqnarray*}
&&E^{P^N} \biggl[ \biggl(\int_0^t \sum_{i<k}
\frac{1}{N^2} \Phi(NM^i,NM^k) \,dL^{ik} \biggr)^2 \biggr]\\
&&\qquad\le
\biggl[1+E^{P^N} \biggl[\int_0^t\sum_{i<k}
\frac{1}{N^3} \Phi(NM^i,NM^k) \,dL^{ik} \biggr] \biggr] \le2 .
\end{eqnarray*}

We take expectations in (\ref{equ36}) and combine with (\ref{equ37}),
(\ref{equ38}) and
(\ref{equ310}) to obtain
\begin{eqnarray*}
&&(1-\Gamma\zeta)E \biggl[\int_0^t \sum_{i<k} M^i M^k\Phi
(NM^i,NM^k)\,dL^{ik} \biggr]
\\
&&\qquad
\le4c\zeta E^{P^N} \biggl[\sum_i M^i[NM^i_t]^{\mathfrak{p}} \biggr] \\
&&\qquad\quad{} +
\Gamma(\zeta) \int_0^t E^{P^N} \biggl[ \biggl(1+\sum_i
M^i[NM^i]^{\mathfrak{p}
} \biggr)
\biggl(1+\sum_i m^i a(Nm^i) \biggr) \biggr] \,ds \\
&&\qquad\quad{} + \Gamma\zeta E^{P^N} \biggl[ \biggl(1+\sum_i
M_t^i[NM_t^i]^{\mathfrak{p}} \biggr)
\int_0^t\sum_{i<k} \frac{1}{N^2} \Phi(NM^i,NM^k)\, dL^{ik} \biggr] \\
&&\qquad\le4c\zeta E^{P^N} \biggl[\sum_i
M^i[NM^i_t]^{\mathfrak{p}} \biggr] \\
&&\qquad\quad{}
+\Gamma(\zeta) E^{P^N} \biggl[1+\sum_i m^i a(Nm^i)^2 \biggr]\int_0^t
E^{P^N} \biggl[1+N\sum_i[M^i]^2 \biggr] \,ds \\
&&\qquad\quad{} + 2 \Gamma\zeta E^{P^N} \biggl[1+N\sum_i [M_t^i]^2
\biggr] .
\end{eqnarray*}
In order to derive this last bound we have used H\"{o}lder's
inequality and the fact that $\mathfrak{p}\le1/2$. From (\ref{equ23}), we
conclude that
%
%
\begin{eqnarray}\label{equ311}
&&(1-\Gamma\zeta) E^{P^N} \biggl[\int_0^t \sum_{i<k} M^i
M^k \Phi(NM^i,NM^k) \,dL^{ik} \biggr] \nonumber\\
&&\qquad\le
\Gamma' \zeta \biggl[1+E^{P^N} \biggl[N\sum_i[M_t^i]^2 \biggr]
\biggr]\\
&&\qquad\quad{} +
\Gamma(\zeta)\int_0^t
E^{P^N} \biggl[1+N\sum_i[M^i]^2 \biggr] \,ds .\nonumber
\end{eqnarray}
Choose $\zeta\le1/(4[\Gamma+\Gamma'])$, where $\Gamma$ and
$\Gamma'$ are the constants appearing in the first line and in front
of the
first term on the right above, respectively.
Combining (\ref{equ311}) with (\ref{equ35}) we get
\begin{eqnarray*}
E^{P^N} \biggl[N\sum_i[M_t^i]^2 \biggr]\le\Gamma\biggl[E^{P^N}
\biggl[1+N\sum_i[m^i]^2
+\int_0^t
\biggl(1+N\sum_i[M^i]^2 \biggr)\,ds \biggr] \biggr].
\end{eqnarray*}
Estimate (\ref{equ31}) now follows from Gronwall's lemma and conditions
(\ref{equ23})
on the initial distribution of masses, and (\ref{equ32}) is
immediate from (\ref{equ311}).
\end{pf*}
\begin{pf*}{Proof of Lemma \protect\ref{lemma2}}
Choose $\zeta>0$ and $\delta>0$ such that
\[
4\zeta[1+\Gamma] [1+C+C(T) ]<\frac{\epsilon
^2}{50}
\quad\mbox{and}\quad
4\delta\Gamma(\zeta) [1+C(T) ] [1+C']\le\frac
{\epsilon
^2}{50} ,
\]
where $\Gamma$ and $\Gamma(\zeta)$ are the constants appearing on the
right of (\ref{equ310}) and (\ref{equ38}), respectively,
$C$ and $C'$ are the constants from assumption (\ref{equ23}), and
$C(T)$ the bound
established in Lemma \ref{lemma1}.
Set the parameter of $g$ equal to a value
of $\zeta$ determined as above.

As in the proof of Theorem \ref{theo1}, we will stop the process at
the finite stopping time $\tau\wedge\tau_{\zeta}\wedge T$, where
$\tau$ is the stopping time defined in the proof of Theorem \ref
{theo1}, and
\[
\tau_{\zeta}=\inf\biggl\{t\ge0, N\sum_i [M_t^i]^2 \ge
\frac{1}{\zeta} \biggr\}.
\]
By Lemma \ref{lemma1} and the choice of $\zeta$, we have
\[
P^N [\tau_{\zeta}\le T ]\le C(T) \zeta\le\frac{\epsilon}{2},
\]
if $\epsilon$ is small enough.
We will thus assume that $P^N$ is
supported on
%
%
\begin{equation}\label{equ312}
\max_i \{M_T^i\} \biggl[1+\sum_i m^i a(Nm^i) \biggr] \le
\frac{1}{N^{1/4}} ,\qquad N\sum_i[M_T^i]^2\le\frac{1}{\zeta} .
\end{equation}
The proof will now follow by estimating the variation of
the terms $\{A_i\}_{1\le i\le5}$ on the right of (\ref{equ36}). Let
$A_i(s,t)=A_i(0,t)-A_i(0,s), 1\le i\le5$.

By Chebyshev's
inequality, (\ref{equ37}), (\ref{equ38}), and the choice of $\zeta,
\delta$, we have
%
%
\begin{eqnarray}\label{equ313}
P^N \biggl[ {\mathop{\mathop{\sup}_{0\le s\le t\le T}}_{t-s\le\delta}}
|A_1(s,t)|>\frac{\epsilon}{5} \biggr]&\le&\frac{\epsilon}{10} ,\\
\label{equ314}
P^N \biggl[ {\mathop{\mathop{\sup}_{0\le s\le t\le T}}_{t-s\le\delta}}
|A_3(s,t)|>\frac{\epsilon}{5} \biggr] &\le&\frac{\epsilon}{10}
\end{eqnarray}
and
%
\begin{equation}\label{equ315}
P^N \biggl[ {\mathop{\mathop{\sup}_{0\le s\le t\le T}}_{t-s\le\delta}}
|A_5(s,t)|>\frac{\epsilon}{5} \biggr] \le\frac{\epsilon}{10} .
\end{equation}

The quadratic variations $Q_2(0,t)$ and $Q_4(0,t)$ of the martingale terms
$A_2$ and $A_4$ satisfy
\begin{eqnarray*}
Q_2(0,t)&\le&\Gamma(\zeta,T) \max_i\{M^i_T\} \biggl( \biggl[\sum_i m^i
a(N m^i) \biggr]
\biggl[1+N\sum_i [M_T^i]^2 \biggr] \nonumber\\
&&\hspace*{85.3pt}{}+
\biggl[ 1+N\sum_i [M^i_T]^2 + \sum_i m^i a(N m^i) \biggr]
\biggr)\nonumber\\
&\le&\frac{\Gamma(\zeta, T)}{N^{1/4}} \biggl[1+\sum_i
m^i a(N m^i) + N\sum_i [M^i_T]^2 \biggr] ,\\
Q_4(0,t)&\le&\Gamma\zeta^2 \biggl[1+N\sum_i[M^i_t]^2 \biggr] \int
_0^t \sum_{i<k}
\biggl[M^i M^k+\frac{1}{N^2} \biggr] \Phi(N M^i,N M^k) \,dL^{ik}\\
&\le&4 \Gamma\zeta\int_0^t \sum_{i<k}
\biggl[M^i M^k+\frac{1}{N^2} \biggr] \Phi(NM^i,NM^k) \,dL^{ik},
\end{eqnarray*}
respectively.
In order to derive these inequalities
we have used the assumptions on the mass sizes, (\ref{equ312}), and the fact
that the diffusion coefficients are decreasing, so that they can be controlled
when the masses $M^i_t$ are large.

Then, by Doob's inequality,
%
%
\begin{eqnarray}\label{equ316}
P^N \biggl[ {\mathop{\mathop{\sup}_{0\le s\le t\le
T}}_{t-s\le\delta}}
|A_2(s,t)|>\frac{\epsilon}{5} \biggr]
&\le&
P^N \biggl[ {2\sup_{0\le t\le T}}
|A_2(0,t)|>\frac{\epsilon}{5} \biggr]\nonumber\\[-8pt]\\[-8pt]
&\le&
\frac{\Gamma(\epsilon,\zeta,T)}{N^{1/4}} [1+C'+C(T) ],\nonumber
\end{eqnarray}
and similarly,
%
%
\begin{eqnarray}\label{equ317}
P^N \biggl[ {\mathop{\mathop{\sup}_{0\le s\le t\le T}}_{t-s\le\delta}}
|A_4(s,t)|>\frac{\epsilon}{5} \biggr]\le
4 \Gamma\zeta [1+C+C(T) ] .
\end{eqnarray}
We pass to the limit $N\rightarrow\infty$ in (\ref{equ316}),
and conclude the proof from the estimates obtained in (\ref{equ313}),
(\ref{equ314}),
(\ref{equ315}), (\ref{equ317}) and the choice of $\zeta$.
\end{pf*}

\section{The hydrodynamic equation}\label{sec4}

We begin with Proposition \ref{prop1}.
\begin{pf*}{Proof of Proposition \protect\ref{prop1}}
Estimate (\ref{equ31}) in Lemma \ref{lemma1} implies that
\[
E^{Q} \biggl[\sup_{t\le T} \langle m\wedge M,\mu_t \rangle\biggr]\le C(T)
\]
uniformly in $M>0$. Then (\ref{equ26}) follows by letting
$M\rightarrow\infty$ and monotone
convergence.

In order to obtain (\ref{equ25}), we will show that
the probability measure $Q$ satisfies
%
%
\begin{equation}\label{equ41}
E^{Q} \biggl[ \mathop{\mathop{\sup}_{\Psi(x,s) \in L^1[\mathbb
{T}\times[0,T]]}}_{
\|\Psi\|_1\le1} \int_0^T \biggl\langle\frac{\omega(m)}{m} \Psi,
\mu_s \biggr\rangle
\,ds \biggr] <\infty.
\end{equation}
Indeed, (\ref{equ41}) implies that with $Q$-probability 1, $\mu_t(dx,dm)=
\upsilon_t(x,dm) \,dx$ for almost every $t\in[0,T]$,
where $\upsilon_t$ satisfies
estimate (\ref{equ25}) in the statement of the proposition. But
$Q$ is supported on $C ([0,T],\mathcaligr{M}_1(\mathbb{T}\times
\mathbb{R}
_+) )$,
hence the result.

We must therefore prove that
\begin{eqnarray*}
&&
E^{Q} \biggl[\sup_{k\in\mathbb{N}}\int_0^T \biggl\langle
\frac{\omega(m)}{m} \Psi^k,\mu_s \biggr\rangle \,ds \biggr]\\
&&\qquad
=\lim_{K\rightarrow\infty}\lim_{\Lambda\rightarrow\infty}
\lim_{N\rightarrow\infty}E^{Q^N} \biggl[\sup_{1\le k\le K}
\int_0^T \biggl\langle \frac{\omega(m)}{m} \Psi^{k,\Lambda},\mu
_s \biggr\rangle \,ds \biggr] <\infty,
\end{eqnarray*}
where we denote
\[
\Psi^{k,\Lambda}(s,x)=\min\{\Psi^k(s,x), \Lambda\},
\]
$\{\Psi^k\}_{k\in\mathbb{N}}$ a dense family in $C
([0,T],\mathbb{T} )
\cap B_1 [L^1([0,T]\times\mathbb{T}) ]$
in the supremum norm, $B_1 [L^1([0,T]\times\mathbb{T}) ]$
the unit ball in $L^1 [[0,T]\times\mathbb{T} ]$. The measures
$\mu_s$ are nonnegative $Q$ a.e., so we may take $\Psi^k\ge0$,
$k \in\mathbb{N}$.

For each $\Psi^{k,\Lambda}$ and $m>0$, let then $u^{k,\Lambda
}(s,x,m)$ be the solution to
\[
\cases{
u^{k,\Lambda}_s+\dfrac{a(m)}{2}u^{k,\Lambda}_{xx}=-\Psi^{k,\Lambda
},\cr
u^{k,\Lambda}(T,\cdot)=0 .}
\]
We have the representation formula
%
%
\begin{equation}\label{equ42}
u^{k,\Lambda}(s,x,m)=\int_s^T\int_{\mathbb{T}} p\bigl(a(m)(u-s),x,z\bigr)
\Psi
^{k,\Lambda}(u,z) \,dz \,du,
\end{equation}
where $p(t,x,z)$ is the Brownian transition density on
$\mathbb{T}$.

It\^{o}'s formula applied to $\frac{\omega(m)}{m}u^{k,\Lambda}$ yields
%
%
\begin{eqnarray}\label{equ43}
&&\int_0^T \sum_i M^i \frac{\omega(NM^i)}{NM^i} \Psi^{k,\Lambda
}(X^i,NM^i) \,ds
\nonumber\\
&&\qquad=
\sum_i \frac{\omega(Nm^i)}{N}
u^{k,\Lambda}(x^i,Nm^i)\nonumber\\[-8pt]\\[-8pt]
&&\qquad\quad{} +\int_0^T \sum_i \frac{\omega(NM^i)}{N} \,\frac
{\partial u^{k,\Lambda}}
{\partial x} \,dX^i\nonumber\\
&&\qquad\quad{} +\int_0^T \sum_{i<j} D_N(X^i,NM^i,NM^j)
(u^{k,\Lambda}) \,dE^{ij},\nonumber
\end{eqnarray}
if $D_N$ now denotes
\begin{eqnarray*}
&&
D_N(x,m,m')(f)\\
&&\qquad
=\frac{1}{N} [\omega(m+m')f(x,m+m')-\omega(m)f(x,m)-\omega(m')f(x,m')].
\end{eqnarray*}
The last term in the expansion (\ref{equ43})
is nonpositive by
formula (\ref{equ42}), the assumption that $\Psi^k \ge0$,
and the subadditivity in $m$ of
$a(m)^{-1/2}\omega(m)$. We thus get
%
%
\begin{eqnarray}\label{equ44}
&&E^{Q^N} \biggl[\sup_{1\le k\le K}
\int_0^T \biggl\langle\frac{\omega(m)}{m} \Psi^{k,\Lambda},\mu
_s \biggr\rangle \,ds \biggr]
\nonumber\\
&&\qquad\le
E^{P^N} \biggl[\sup_{1\le k\le K}
\sum_i \frac{\omega(Nm^i)}{N} u^{k,\Lambda}(x^i,Nm^i) \biggr]
\\
&&\qquad\quad{} +
E^{P^N} \biggl[\sup_{1\le k\le K} \biggl|
\int_0^T \sum_i \frac{\omega(NM^i)}{N} \,\frac{\partial
u^{k,\Lambda}}
{\partial x} \,dX^i \biggr| \biggr].\nonumber
\end{eqnarray}
The second term on the right-hand side above can be easily bounded by
replacing the supremum by a sum over $1\le k \le K$ and computing the
quadratic variation of each of the resulting orthogonal martingale
terms. The sum of these quadratic variations vanishes in the limit
$N\rightarrow\infty$; in order to see this, it suffices to replace
$u^{k,\Lambda}$ by its representation (\ref{equ42}) and then apply
assumptions (\ref{equ21b}), (\ref{equ23}) and estimate~(\ref{equ31}).

Finally, the hypothesis on $P^N_0$, $\nu$, $\nu^*$ and the fact that
$\|\Psi^{k,\Lambda}\|_1\le\|\Psi^k\|_1 \le1$ imply that
\begin{eqnarray*}
&&\lim_{N\rightarrow\infty} E^{P^N} \biggl[ \sup_{1\le k\le K}
\sum_i \frac{\omega(Nm^i)}{N} u^{k,\Lambda}(x^i,Nm^i) \biggr]\\
&&\qquad=\sup_{1\le k \le K}
\biggl\langle
\frac{\omega(m)}{m} u^{k,\Lambda},\nu \biggr\rangle\le
\sup_{1\le k\le K} \biggl\langle
\frac{\omega(m)}{m},\nu^* \biggr\rangle\|\Psi^{k,\Lambda}\|_1
=\Gamma(\nu^*)<\infty
\end{eqnarray*}
holds uniformly in $\Lambda, K$.

We pass to the limit $\Lambda\rightarrow\infty$ and
then $K \rightarrow\infty$ in (\ref{equ44})
to obtain (\ref{equ41}).
\end{pf*}
\begin{pf*}{Proof of Theorem \protect\ref{theo2}}
It will be enough to consider $f \in C_b^2(\mathbb{T}\times\mathbb{R}_+)$
compactly supported and Lipschitz in $m$, and then use bounded
convergence to obtain (\ref{equ27})
for a general $f\in C_b^2(\mathbb{T}\times\mathbb{R}_+)$. We need to
analyze the difference
\[
Z_f(t)= \sum_i M^i_t f(X^i_t,NM^i_t)-\sum_i m^i f(x^i,Nm^i) .
\]
We start by writing the semimartingale $Z_f$ as
\[
Z_f(t)= \mathcaligr{H}_f(t)+A_f(t),
\]
where $\mathcaligr{H}_f$ is the martingale obtained by adding the fluctuation
terms arising from the free particle dynamics and the stochastic
coagulation phenomena. These can be proved negligible by applying
Doob's
inequality, the integrability assumptions on $a(m)$ stated in Section
\ref{sec1}, and Lemma \ref{lemma1}. The term $A_f$ is given by
\begin{eqnarray*}
A_f(t) &=& \int_0^t \sum_i M^i_s \frac{a(NM^i)}{2} \frac{\partial^2
f}{\partial x^2} (X^i_s,NM^i_s) \,ds \\
&&{}
+ \int_0^t \sum_{i<j}
\bigl[M^i \bigl( f\bigl(X^i,N(M^i+M^j)\bigr)-f(X^i,NM^i) \bigr)\\
&&\hspace*{44pt}{} + M^j \bigl(f\bigl(X^j,N(M^i+M^j)\bigr)-f(X^j,NM^j) \bigr) \bigr]
\\
&&\hspace*{39.1pt}{}\times
\frac{\Phi(NM^i,NM^j)}{N} \,dL^{ij}.
\end{eqnarray*}

The first term of $A_f$ will clearly have the limit
%
%
\begin{equation}\label{equ45}
\int_0^t \int_{\mathbb{T}\times\mathbb{R}_+}
\biggl\langle\frac{a(m)}{2} \frac{\partial^2
f}{\partial x^2},\mu_s \biggr\rangle \,ds.
\end{equation}
We can guess the limit of the second term from the occupation times
formula; we should recover the second term in the hydrodynamic equation
(\ref{equ27}). In order to obtain this expression, we replace $dL^{ij}$ in
the second term of $A_f$ by $[a(m)+a(m')] V_{\epsilon} (X^i-X^j) \,ds$,
where $V_{\epsilon}(x)$
approximates the Dirac $\delta$-function at the origin as $\epsilon
\rightarrow0$. The new
integral will converge
weakly to
\begin{eqnarray*}
&&\frac{1}{2} \int_0^t \int_{\mathbb{T}\times\mathbb{T}}
\int_{\mathbb{R}_+ \times\mathbb{R}_+}
\biggl[\frac{f(x,m+m')-f(x,m)}{m'}+\frac{f(y,m+m')-f(y,m')}{m}
\biggr]\\
&&\hspace*{60.7pt}\qquad{}
\times\kappa(m,m') V_{\epsilon}(x-y) \mu_s(dx,dm) \mu_s(dy,dm')
\,ds
\end{eqnarray*}
as $N\rightarrow\infty$.

We will justify this exchange by showing that there exists a sequence
of measurable
sets $\mathcaligr{C}_{N,\epsilon,T}$ with
\[
\lim_{\epsilon\rightarrow0}\limsup_{N\rightarrow\infty}
P^N [\mathcaligr{C}_{N,\epsilon,T} ]=0,
\]
such that
$\Upsilon_{N,\epsilon,f}(t)$ given by
\begin{eqnarray*}
&&\int_0^t \sum_{i<j}
\bigl[M^i \bigl( f\bigl(X^i,N(M^i+M^j)\bigr)-f(X^i,NM^i) \bigr)\\
&&\qquad\hspace*{8.5pt}{}
+M^j \bigl(f\bigl(X^j,N(M^i+M^j)\bigr)-f(X^j,NM^j) \bigr) \bigr]
\frac{\Phi(NM^i,NM^j)}{N}\\
&&\qquad\hspace*{4pt}{}\times\bigl(dL^{ij}- [a(NM^i)+a(NM^j) ]
V_{\epsilon
}(X^i_s-X^j_s) \,ds \bigr)
\end{eqnarray*}
satisfies
%
%
\begin{equation}\label{equ46}
\lim_{\epsilon\rightarrow0} \limsup_{N\rightarrow\infty}
E^{P^N} \biggl[{\sup_{0\le t\le T}} | \Upsilon_{N,\epsilon
,f}(t) |
1_{\mathcaligr{C}^c_{N,\epsilon,T}} \biggr]=0 .
\end{equation}
Here $\mathcaligr{C}^c$ denotes the complement of the set $\mathcaligr{C}$.

Suppose that (\ref{equ46}) holds, and let $\delta>0, l \in\mathbb{N}$. Define
\begin{eqnarray*}
\kappa^l(m,m') &=& \kappa(m,m') 1_{\{l^{-1}\le m \le l\}}
1_{\{l^{-1}\le m' \le l\}},\\
a^l(m) &=& a(m) 1_{\{l^{-1}\le m \le l\}}
\end{eqnarray*}
and
\begin{eqnarray*}
\mathcaligr{F}_{l,\epsilon,\delta} &=& \biggl\{\sup_{0\le t\le T}
\biggl| \langle\mu_t,f \rangle
-\langle\mu_0,f \rangle-\int_0^t \biggl\langle\frac{1}{2} a^l(m)
\,\frac{\partial^2 f}{\partial x^2},\mu_s \biggr\rangle \,ds \\
&&\hspace*{33.7pt}{} -\int_0^t
\int_{(\mathbb{T}\times\mathbb{R}_+)^2} \frac
{[f(x,m+m')-f(x,m)]}{m'}
V_{\epsilon}(x-y)\\
&&\hspace*{100.02pt}{} \times
\kappa^l(m,m') \mu_s(dx,dm) \mu_s(dy,dm') \,ds
\biggr|\le\delta \biggr\}.
\end{eqnarray*}
Then $\mathcaligr{F}_{l,\epsilon,\delta}$ is closed in
$C([0,T],\mathcaligr{M}_1(\mathbb{T}\times\mathbb{R}_+))$
with the Skorokod topology.

By Proposition \ref{prop1}, $Q$ almost everywhere,
\[
\lim_{l\rightarrow\infty} \int_0^t \biggl\langle\frac{1}{2}
a^l(m)
\,\frac{\partial^2 f}{\partial x^2},\mu_s \biggr\rangle \,ds
=\int_0^t \biggl\langle\frac{1}{2} a(m)
\,\frac{\partial^2 f}{\partial x^2},\mu_s \biggr\rangle \,ds
\]
and
\begin{eqnarray*}
&&\lim_{l \rightarrow\infty}
\int_0^t \biggl\langle\frac{[f(x,m+m')-f(x,m)]}{m'}
V_{\epsilon}(x-y) \kappa^l(m,m'),\mu_s \otimes\mu_s \biggr\rangle
\,ds\\
&&\qquad
=\int_0^t \biggl\langle\frac{[f(x,m+m')-f(x,m)]}{m'}
V_{\epsilon}(x-y) \kappa(m,m'),\mu_s \otimes\mu_s \biggr\rangle \,ds
\end{eqnarray*}
for all $t\in[0,T]$. These imply
%
%
\begin{eqnarray}\label{equ47}
\limsup_{l\uparrow\infty}
\limsup_{N\uparrow\infty}
Q^N [\mathcaligr{F}_{l,\epsilon,\delta} ]
&\le& \limsup_{l\uparrow\infty}
Q [ \mathcaligr{F}_{l,\epsilon,\delta} ]
\nonumber\\[-8pt]\\[-8pt]
&\le&
Q [\mathcaligr{F}_{\infty,\epsilon,2\delta} ].\nonumber
\end{eqnarray}

Now, we know that $\mu$
disintegrates
as $\mu_s(dx,dm)=\upsilon_s(x,dm) \,dx$, $s \in[0,T]$, $Q$ a.e.
Letting $\epsilon\rightarrow
0$ in (\ref{equ47}), by Lebesgue's differentiation
theorem, dominated convergence and (\ref{equ46}), we have
\begin{eqnarray*}
1&=&\lim_{\epsilon\to0} \limsup_{l\uparrow\infty}
\limsup_{N\uparrow\infty}
Q^N [\mathcaligr{F}_{l,\epsilon,\delta} ] \le\lim
_{\epsilon\to0}
Q [\mathcaligr{F}_{\infty,\epsilon,2\delta} ]
\\
&=&
Q [\mathcaligr{F}_{\infty,2\delta} ],
\end{eqnarray*}
if $\mathcaligr{F}_{\infty,2\delta}$ is obtained replacing
$V_{\epsilon}(x-y)$ in the definition of $\mathcaligr{F}_{\infty,
\epsilon, 2\delta}$
by the Dirac function evaluated at $x-y$.
Since $\delta>0$ is arbitrary, this implies that (\ref{equ27})
holds with $Q$-probability 1.

It remains to prove (\ref{equ46}). For each $\epsilon>0$, let the
approximation of
the Dirac $\delta$
function $V_{\epsilon}$ be such that there exists
a function $u_{\epsilon}$ in
$C^2(\mathbb{T}-\{0\}) \cap C^1(\mathbb{T})$ with support
contained in
$(-1/2,1/2)$, so that
\begin{eqnarray*}
\|u_{\epsilon}\|_{\infty}&\le&\gamma_1 \epsilon,\qquad u'_{\epsilon}(0) = \frac{1}{2},
\\
\|u_{\epsilon}'\|_{\infty}&\le& \gamma_2,\qquad \lim_{\epsilon\rightarrow
0}u_{\epsilon}'(x)=0,\qquad x\neq0,
\end{eqnarray*}
and
\[
u_{\epsilon}''(x)=W_{\epsilon}(x),\qquad x\neq
0;\qquad u_{\epsilon}''(x)<0 \qquad\mbox{if } 0<x<\epsilon,
\]
where $W_{\epsilon}$
is a real valued function such that
$|W_{\epsilon}|=V_{\epsilon}$, and $\gamma_1$ and $\gamma_2$ are positive
constants independent of $\epsilon$.

Consider the finite stopping time
\[
\tau_{\epsilon}=\inf\biggl\{t\le T, N\sum_i
[M^i_t]^2 \ge\frac{1}{\sqrt{\epsilon}} \biggr\}
\]
and define
\[
\mathcaligr{C}_{N,\epsilon,T}= \{\tau_{\epsilon}\le T \} .
\]
Lemma \ref{lemma1} then implies
\[
P^N [\mathcaligr{C}_{N,\epsilon,T} ]\le C(T) \sqrt
{\epsilon}
\]
and clearly $\lim_{\epsilon\rightarrow0}\limsup_{N\rightarrow
\infty}
P^N[\mathcaligr{C}_{N,\epsilon,T}]=0$.

Let
\begin{eqnarray*}
G(x,y,m,m') &=& G_{N,f,\epsilon}(x,y,m,m')\\
&=& \bigl[ m
[f(x,m+m')-f(x,m)
]\\
&&\hspace*{2pt}{} +m' [f(y,m+m')-f(y,m') ] \bigr]
u_{\epsilon}(|x-y|) \frac{\Phi(m,m')}{N^2} .
\end{eqnarray*}
The proof follows the usual pattern after this point: we will stop the
process upon achieving the stopping time $\tau_{\epsilon}$, so that
we may
assume that $P^N$ is supported on
\[
N\sum_i [M^i_T]^2 \le\frac{1}{\sqrt{\epsilon}}.
\]
We will then
apply It\^{o}--Tanaka's formula to
\[
\sum_{i<j} G(X^i_t,X^j_t,NM^i_t,NM^j_t)
\]
in order to recover $\Upsilon_{N,\epsilon,f}$ from the nondifferentiability
of $u_{\epsilon}(|x|)$ at the
origin and the particular choice of $u_{\epsilon}$.
One then has to check that the remaining terms of the expansion
converge to $0$ uniformly on $[0,T]$, when
taking $N\rightarrow\infty$ and $\epsilon\rightarrow0$, in that
order.

We have
\begin{eqnarray*}
&& \sum_{i<j} G(X^i_t,X^j_t,NM^i_t,NM^j_t)\\
&&\qquad =
\sum_{i<j} G(x^i,x^j,Nm^i,Nm^j)+\frac{1}{2} \Upsilon_{N,\epsilon,f}(t)
+\mathcaligr{H}_{N,\epsilon,f}(t)+
A_{N,\epsilon,f}(t),
\end{eqnarray*}
where $\mathcaligr{H}_{N,\epsilon,f}(t)$ is a $P^N$-martingale and
$A_{N,\epsilon
,f}(t)$ is a
continuous, bounded variation process. We start with the former:
\begin{eqnarray*}
\mathcaligr{H}_{N,\epsilon,f}
&=&\sum_{i} \int_0^t \sum_{j\neq i} \,\frac{\partial
G}{\partial X^i} (X^i_s,X^j_s, NM^i_s,NM^j_s)
\,dX^i\\
&&{} +
\sum_{i<k} \int_0^t
\biggl( \biggl[ \sum_{j\neq i,k} G(N[M^i+M^k],NM^j)\\
&&\hspace*{51.72pt}{} - G(NM^i,NM^j)-G(NM^k,NM^j) \biggr]\\
&&\hspace*{142.25pt}{}-G(NM^i,NM^k) \biggr)\\
&&\hspace*{40pt}{} \times\biggl[dE^{ik}
-\frac{\Phi(NM^i,NM^k)}{N}\,dL^{ik} \biggr].
\end{eqnarray*}
The quadratic variation of the first term on the right above, the
Brownian martingale, can be easily seen to vanish when $N\rightarrow
\infty$. We proceed to show how to bound
one
term in the quadratic variation $Q_c$ of the coagulation martingale.
Consider then
\begin{eqnarray*}
Q_c^1(t)
&=&\int_0^t\sum_{i<k} [M^i]^2 \biggl\{ \sum_{j\neq i,k} \bigl[f
\bigl(X^i,N(M^i+M^k+M^j) \bigr)
\\
&&\hspace*{84.24pt}{}
-f \bigl(X^i,N(M^i+M^k) \bigr)\\
&&\hspace*{84.24pt}{}-f \bigl(X^i,N(M^i+M^j) \bigr)+f(X^i,NM^i) \bigr]\\
&&\hspace*{80pt}{}\times u_{\epsilon}(|X^i-X^j|) \frac{\Phi
(N(M^i+M^k),NM^j )}{N} \biggr\}^2\\
&&\hspace*{53.36pt}{}\times
\frac{\Phi(NM^i,NM^k)}{N} \,dL^{ik}.
\end{eqnarray*}
We now use that $f$ has compact support. Let $L\ge0$ be such that
$f(x,m)=0$ whenever $|m|>L$. The expression between brackets in $Q^1_c$
will thus vanish whenever $NM^i>L$, so that we may bound $Q_c^1(T)$ by
\begin{eqnarray*}
&&\Gamma(f) \epsilon^2\int_0^T \sum_{i<k} \frac{1}{N^2}
\biggl[ [N(M^i+M^k) ]^{2\mathfrak{p}}+\sum_j
\frac{1}{N}[NM^j]^{2\mathfrak{p}} \biggr]\\[-0.5pt]
&&\hspace*{66.5pt}{}\times\frac{\Phi
(NM^i,NM^k)}{N}\,dL^{ik}\\[-0.5pt]
&&\qquad\le\Gamma(f) \epsilon^2 \int_0^T \sum_{i<k}
\frac{1}{N^2} \Phi(NM^i,Nm^k)\,dL^{ik}
\end{eqnarray*}
by the assumption that $\mathfrak{p}\le1/2$. Now,
\begin{eqnarray*}
&&E^{P^N} \biggl[\int_0^T \sum_{i<k}
\frac{1}{N^2} \Phi(NM^i,NM^k) \,dL^{ik} \biggr]\\[-0.5pt]
&&\qquad=E^{P^N} \biggl[\sum_i
\frac{1}{N}1_{\{m^i>0\}}-\sum_i \frac{1}{N}1_{\{M^i_T>0\}} \biggr]
\le1 ,
\end{eqnarray*}
from where it follows that $\lim_{\epsilon\rightarrow0}\lim
_{N\rightarrow
\infty} E^{P^N}[Q^1_c]=0$. Similar considerations prove that the
expectation of the rest
of the terms in $Q_c$ vanish in the limit $N\rightarrow\infty$,
$\epsilon\rightarrow0$.
Doob's inequality then implies
%
%
\begin{equation}\label{equ48}
\lim_{\epsilon\rightarrow0}\limsup_{N\rightarrow\infty}
E^{P^N} \biggl[{\sup_{0\le t\le T}} |\mathcaligr{H}_{N,\epsilon,f}(t)|
1_{C^c_{N,\epsilon,T}} \biggr]=0 .
\end{equation}

The Lipschitz property of $f$ yields
%
%
\begin{eqnarray}\label{equ49}
&& E^{P^N} \biggl[
\biggl|\sum_{i<j} G(X_t^i,X_t^j,M^i_t,M_t^j)-\sum_{i<j}
G(x^i,x^j,m^i,m^j) \biggr| \biggr]\nonumber\\[-8pt]\\[-8pt]
&&\qquad
\le\Gamma(f) \epsilon E^{P^N} \biggl[1+N\sum_i
[M_t^i]^2 \biggr]\le\Gamma(f) \epsilon[1+C(T)].\nonumber
\end{eqnarray}

The process $A_{N,\epsilon,f}(t)$ equals
\begin{eqnarray*}
A_{N,\epsilon,f}(t) &=& \frac{1}{2}\int_0^t \sum_{i\neq j} M^i
\biggl[\frac{\partial^2 f}{\partial x^2}\bigl(X^i,N(M^i+M^j)\bigr) -
\frac{\partial^2 f}{\partial x^2}(X^i,NM^i) \biggr]
\nonumber\\[-0.5pt]
&&\hspace*{52.43pt}{} \times u_{\epsilon}(|X^i-X^j|) \frac{\Phi
(NM^i,NM^j)}{N}
a(NM^i) \,ds
\nonumber\\[-0.5pt]
&&{} + \frac{1}{2} \int_0^t \sum_{i\neq j} M^i
\biggl[\frac{\partial
f}{\partial x} \bigl(X^i,N(M^i+M^j)\bigr)-\frac{\partial f}{\partial
x}(X^i,NM^i)
\biggr]\\[-0.5pt]
&&\hspace*{65.3pt}{} \times\operatorname{sign}(X^i-X^j)
u_{\epsilon}'(|X^i-X^j|) \\[-0.5pt]
&&\hspace*{65.3pt}{} \times\frac{\Phi(NM^i,NM^j)}{N} a(NM^i) \,ds \\
&&{} +\int_0^t \biggl(\sum_{i<k} \sum_{j\neq i, k}
\bigl[ G\bigl(X^i,X^j,N(M^i+M^k),NM^j\bigr)\\
&&\hspace*{71.3pt}{} -
G(X^i,X^j,NM^i,NM^j)\\
&&\hspace*{94pt}{} - G(X^k,X^j,NM^k,NM^j) \bigr]\\
&&\hspace*{101.4pt}{} - G(X^i,X^k,NM^i,NM^k) \biggr)\\
&&\hspace*{24.1pt}{}\times
\frac{\Phi(NM^i,NM^k)}{N} \, dL^{ik}\\
&=& I^1_{N,\epsilon,f}(t)+I^2_{N,\epsilon
,f}(t)+I^3_{N,\epsilon,f}(t),
\end{eqnarray*}
where $\operatorname{sign}(x)$ takes values $1$ or $-1$ according to whether
$x> 0$ or $x\le0$. It is easy to see that
%
%
\begin{equation}\label{equ410}
\lim_{\epsilon\rightarrow0}\limsup_{N\rightarrow\infty}
E^{P^N} \biggl[{\sup_{0\le t\le T}}|I^1_{N,\epsilon,f}(t)| \biggr]=0 .
\end{equation}
We can then bound
\begin{eqnarray*}
|I_{N,\epsilon,f}^2(t) |&\le&\Gamma(f) \int_0^t \sum
_{i<j} M^i
|u_{\epsilon}'|(|X^i-X^j|) \frac{\Phi(NM^i,NM^j)}{N} a(NM^i) \,ds\\
&\le&\Gamma(f) \biggl[ \sum_i M^i_t [NM^i_t]^{\mathfrak{p}}+\sum_i m^i
a(Nm^i) \biggr] .
\end{eqnarray*}
Since $\lim_{\epsilon\rightarrow0} u_{\epsilon}'(x)=0$ for all
$x\neq0$
in $\mathbb{T}$, dominated convergence implies that
%
%
\begin{equation}\label{equ411}
\lim_{\epsilon\rightarrow0} \limsup_{N\rightarrow\infty}
E^{P^N} \biggl[{\sup_{0\le t\le T} }|I_{N,\epsilon,f}^2(t)| \biggr]=0
\end{equation}
as well. Finally, the fact that $f$ is Lipschitz and (\ref{equ21a}) yield
\begin{eqnarray*}
|I_{N,\epsilon,f}^3(t) | &\le& \Gamma(f,\Phi) \epsilon
\biggl[\sum_i M_t^i(NM_t^i)^{\mathfrak{p}} \biggr]\\
&&{}
\times\int_0^t\sum_{i<k}
\biggl[M^iM^k+\frac{1}{N^2} \biggr] \Phi(NM^i,NM^k) \, dL^{ik}
\end{eqnarray*}
so that
%
%
\begin{equation}\label{equ412}
\lim_{\epsilon\rightarrow0} \limsup_{N\rightarrow\infty}
E^{P^N} \biggl[{\sup_{0\le t\le T}} |I_{N,\epsilon,f}^3(t)|
1_{C^c_{N,\epsilon
,T}} \biggr]=0.
\end{equation}

The limit (\ref{equ46}) is immediate from estimates (\ref{equ48}), (\ref
{equ49}), (\ref{equ410}),
(\ref{equ411}) and (\ref{equ412}), and the result follows.
\end{pf*}

\section{Uniqueness of the solution}\label{sec5}

In this section we seek to establish the uniqueness in
an appropriately defined class
of the solution to the hydrodynamic equa\-tion~(\ref{equ27}), which in
differentiated form can be written as
%
%
\begin{equation}\label{equ51}
\dot{\mu}_t =\frac{1}{2} a(m) \,\frac{\partial^2}{\partial x^2}\,\mu
_t +K(\mu_t), \qquad \mu_0=\nu.
\end{equation}
Here\vspace*{1pt} $\frac{\partial^2 }{\partial x^2}$ denotes the
partial derivative of second order with respect to the space variable,
interpreted in the weak sense, and $K$ is the coagulation kernel given
by
\begin{eqnarray*}
\langle f,K(\mu) \rangle
&=&
\int_{\mathbb{R}_+\times\mathbb{R}_+}
[f(m+m')-f(m) ] \frac{\kappa(m,m')}{m'}
\mu(dm) \mu(dm')
\\
&=&
\frac{1}{2} \int_{\mathbb{R}_+^2}
[(m+m') f(m+m')-m f(m)-m'f(m') ]\\
&&\hspace*{24.3pt}{}\times \frac{\kappa(m,m')}{m m'}
\mu(dm)\mu(dm'),
\end{eqnarray*}
if $f$ is a bounded test function,
$\mu\in\mathcaligr{M}_f(\mathbb{R}_+)$ a finite measure such that
$\kappa(m,\break m')/m'\in L^1(d\mu\times d\mu)$.

We will work under the assumption that there exists a pair of functions
$\varpi$ and $\omega$ bounded on each compact subset of $(0,\infty )$,
such that $\varpi\omega^{-1}$ is bounded, $\omega$, $a^{-1/2} \omega$
and $a^{-1/2} \omega\varpi$ are subadditive, and such that the
following inequalities hold:
%
%
\begin{eqnarray}
\label{equ52}
\kappa(m,m') &\le& \omega(m) \omega(m'),\\
\label{equ53}
\kappa(m,m') &\le& \omega(m)\varpi(m')+\varpi(m)\omega(m').
\end{eqnarray}
We will also require that $a$ be nonincreasing and
the initial measure $\nu$ satisfy
%
%
\begin{equation}\label{equ54}
\nu(dx,dm) \le\nu^*(dm) \,dx,\qquad \nu^*\mbox{ such that }
\biggl\langle \frac{\omega^2}{m},\nu^* \biggr\rangle<\infty.
\end{equation}
Let us define the class $\mathcaligr{B}(\omega)$ by
\begin{eqnarray*}
\mathcaligr{B}(\omega) &=& \biggl\{\eta\in C \bigl(\mathbb
{R}_+,\mathcaligr
{M}_f(\mathbb{T}\times\mathbb{R}_+ ) \bigr)\dvtx
\eta_t(dx,dm)=\upsilon_t(x,dm) \,dx, \\
&&\hspace*{6.5pt}\sup_{t\ge0}
\biggl\| \biggl\langle\!\biggl\langle\frac{\omega(m)}{m},
\upsilon_t \biggr\rangle\!\biggr\rangle\biggr\|_{\infty}
<\infty
\mbox{ and } \sup_{t\le T} \langle m, \mu_t \rangle
<\infty\ \forall
T>0
\biggr\}.
\end{eqnarray*}

The choices
\[
\omega(m)=[1+c'+a(1)] [m^{\mathfrak{p}}+a(m)+1 ]
\quad\mbox{and}\quad \varpi(m)=a(m)+1
\]
satisfy conditions (\ref{equ52}) and (\ref{equ53}) in the particular
situation treated in
Sections \ref{sec3} and \ref{sec4}. In this
case
the coagulation kernel is given by
\[
\kappa(m,m')=\Phi(m,m') [a(m)+a(m')]
\]
with symmetric rate $\Phi$ verifying
\[
\Phi(m,m')\le c (m^{\mathfrak{p}}+m'^{\mathfrak
{p}})1_{[m>0,m'>0]},\qquad 0\le\mathfrak{p}
\le\tfrac{1}{2},
\]
and diffusivity $a(m)$ such that $a^{-1/2}\omega$ is subadditive. Then
Theorem \ref{theo2} yields an
existence result in $\mathcaligr{D}(\omega)$. The uniqueness of this
solution in the larger class
$\mathcaligr{B}(\omega)$ follows from Theorem \ref{theo3} below.

Here is the main result of this section:
\begin{theorem}\label{theo3}
Assume conditions (\ref{equ52}), (\ref{equ53}) and (\ref{equ54}) on the
coagulation rate
$\kappa$ and the initial measure $\nu$. Then for any $T>0$, (\ref
{equ51}) has at
most one solution $\{\mu_t\}_{0\le t\le T}$ in $\mathcaligr{B}(\omega)$.
\end{theorem}

We state the theorem on $\mathbb{T}$ to avoid
introducing more terminology; in fact, the proof holds in $\mathbb{R}^d$
for a general diffusion model with coefficients given by $a$,
undergoing coagulation
at a rate determined by $\kappa$. In that case we require that
(\ref{equ52}), (\ref{equ53}) and (\ref{equ54})
be satisfied by a pair of maps $\omega, \varpi$ such that
$\varpi^{-1}\omega$ is bounded, as before, and $a^{-d/2}\omega$ and
$a^{-d/2}\omega\varpi$
are both subadditive.

The result will follow from considering an approximating system of
equations to the coagulation equation, for which existence and
uniqueness can be easily derived. The method is a simplification of a
technique developed by Norris in \cite{NorrisBcoag}, we have thus
tried to adhere to his notation whenever possible. We are able to make a
significant shortcut in the proof
due to the assumption that we already have got one
solution in $\mathcaligr{B}(\omega)$, this yields a crucial
monotonicity property in the approximating scheme as a direct byproduct
of its construction (compare Lemma 5.1 below with Lemmas 5.5 and 5.6 in
\cite{NorrisBcoag}).

Before we can proceed to prove the theorem we need to
introduce some definitions.
Given $s>0$ and $\mu\in\mathcaligr{M}_f(\mathbb{T}\times\mathbb{R}_+)$,
we will denote $P_s\mu\in\mathcaligr{M}_f (\mathbb{T}\times\mathbb
{R}_+)$ the
measure given by
\[
\langle f,P_s\mu\rangle=\int f(z,m) p(a(m)s,x,z) \mu(dx,dm) \,dz,
\]
where, as before, $p(t,x,z)$ is the Brownian transition density
on $\mathbb{T}$. Hereafter, we will use the abridged notation $p_t^{x,z}(m)$
to denote $p(a(m)t,x,z)$.
We also introduce
the kernels $K^{+}$ and $K^{-}$ on $\mathcaligr{M}_f(\mathbb{R}_+)$
defined as
\begin{eqnarray*}
K^+(\mu)(dm) &=& \int_{m'+m''=m}
\frac{\kappa(m',m'')}{m''} \mu(dm') \mu(dm''),\\
K^-(\mu)(dm) &=& \int_{\mathbb{R}_+} \frac{\kappa(m,m')}{m'} \mu
(dm) \mu(dm').
\end{eqnarray*}
With this notation, we now show that (\ref{equ51}) is equivalent for
$\mu\in
\mathcaligr{B}(\omega)$ to the integral equation
%
%
\begin{equation}\label{equ55}
\mu_t=P_t\nu+\int_0^t P_{t-r}K^+(\mu_r)-P_{t-r}K^-(\mu_r) \,dr.
\end{equation}
By integrating against a test function and differentiating in time, it
follows that any solution to (\ref{equ55}) satisfies (\ref{equ51}).
Conversely, let
$\mu_t$ be a solution to (\ref{equ51}), and set $\tilde{\mu}_t$ equal
to the
right-hand side of (\ref{equ55}).
Then $\mu-\tilde{\mu}_t$ verifies
\[
\frac{d}{dt}
(\tilde{\mu}_t-\mu_t)=\frac{1}{2} a(m) \,\frac{\partial^2
}{\partial x^2}(\tilde{\mu}_t-\mu_t),\qquad
\tilde{\mu}_0-\mu_0=0,
\]
and we conclude that $\tilde{\mu}_t=\mu_t$. In particular $\mu_t$
is a
solution to (\ref{equ55}).

Given a bounded map
$c_s(x)\dvtx[0,T]\times\mathbb{T}\longrightarrow\mathbb{R}$,
let
$\tilde{p}_m(c)$ be the propagator associated with the operator
$\frac{1}{2} a(m) \frac{\partial^2}{\partial x^2}-c_s(\cdot)$ on
$\mathbb{T}$, and consider the kernel
\[
\tilde{P}_{ts}(c)\mu(x,dm)=\int_{\mathbb{T}}
\mu(dz,dm) \tilde{p}_m(c)(s,z;t,x) \,dz,
\]
$\mu$ in
$\mathcaligr{M}_f(\mathbb{T}\times\mathbb{R}_+)$.

Let $E_n=[1/n,n]$ for $n\in\mathbb{N}$, and define $K_n^+$, $K_n^-$
by analogy with $\mathop{\mathop{K}^{+}}^{-}\dvtx\break
K_n^+(\mu)(dm)=K^+(\mu)(dm) 1_{\{m\in E_n\}}$,
\[
K_n^-(\mu)(dm)=\int_{\mathbb{R}_+} \frac{\kappa(m,m')}{m'}
1_{\{m+m'\in E_n\}}
\mu(dm) \mu(dm')
\]
and $K_n=K_n^+ - K_n^-$. Define $\nu^n$ as the restriction of $\nu$
to $E_n$,
$\nu^n(dm)=\nu(dm)\times 1_{\{m\in E_n\}}$.

In order to keep the number of definitions to a minimum,
given a test function $g$, we will
use $\llangle g,\mu\rrangle$ to denote the
integral $\llangle g,\nu\rrangle$ in the case that the decomposition
$\mu(dx,dm)=\nu(x,dm) \,dx$ holds, where we recall that
$\llangle g,\nu\rrangle$ stands for the single integral $\int g(m) \nu(x,dm)$.

The proof of Theorem \ref{theo3} will follow from the following lemma.
\begin{lemma}\label{lemma3}
Let $\{\mu_s\}_{0 \le s \le T} \in\mathcaligr{B}(\omega)$ be a solution
to (\ref{equ51}) with initial value
$\nu$, and assume that conditions (\ref{equ52}), (\ref{equ53}) and (\ref
{equ54}) hold.
Then, for each $n\in\mathbb{N}$,
there exists a unique kernel $(\tilde{\mu}^n_t)_{0\le t \le T}$
in $\mathcaligr{B}(\omega)$, such that
%
%
\begin{equation}\label{equ56}
\tilde{\mu}^n_t=\tilde{P}_t(c^n)\nu^n+ \int_0^t \tilde
{P}_{ts}(c^n) [
K_n^+(\tilde{\mu}_s^n) +\delta^n_s
\tilde{\mu}^n_s ] \,ds,
\end{equation}
where
\begin{eqnarray*}
\delta^n_t(x,m) &=& \biggl\langle\!\biggl\langle\frac{\omega(m) \omega(m')-\kappa
(m,m')}{m'},\tilde{\mu}^n_t\biggr\rangle\!\biggr\rangle, \\
c^n_t(x,m)
&=&
\biggl\langle\!\biggl\langle\omega(m)\frac{\omega(m')}{m'},P_t\nu\biggr\rangle\!\biggr\rangle
+\int_0^t\biggl\langle\!\biggl\langle\omega(m)\frac{\omega(m')}{m'},P_{t-s} [K_n(\tilde
{\mu}
_s^n) ]
\biggr\rangle\!\biggr\rangle \,ds.
\end{eqnarray*}
Moreover, $\tilde{\mu}_s^n$ satisfies
\[
\tilde{\mu}_s^n \le\tilde{\mu}_s^{n+1}
\le\mu_s \qquad\mbox{for all } s\in[0,T].
\]
\end{lemma}
\begin{pf}
The method of the proof is classical. We will take advantage of the fact
that all relevant quantities are bounded in $E_n$ to define
a Picard iteration procedure,
and then prove by a contraction mapping argument
that the scheme converges to a solution, which is finally
shown to be unique.

Define
\begin{eqnarray*}
\tilde{\mu}_{t,0}^n &=& \mu_t 1_{\{m\in E_n\}}\le\tilde{\mu
}_{t,0}^{n+1}\le\mu_t,\\
c_{t,0}^n &=& \biggl\langle\!\biggl\langle\omega(m)\frac{\omega(m')}{m'},P_t \nu\biggr\rangle\!\biggr\rangle
+\int_0^t \biggl\langle\!\biggl\langle\omega(m)\frac{\omega(m')}{m'},
P_{t-s} [K_n(\tilde{\mu}_{s,0}^n) ] \biggr\rangle\!\biggr\rangle \,ds.
\end{eqnarray*}
Condition (\ref{equ54}) on $\nu$ and the fact that
$\tilde{\mu}^n_{\cdot,0}$ is supported on $E_n$ imply that
$c^n_{t,0}(x,m)$ is well defined and bounded.

Let
\[
c_t(x,m)=\biggl\langle\!\biggl\langle\omega(m)\frac{\omega(m')}{m'}, \mu_t \biggr\rangle\!\biggr\rangle.
\]
We claim that $0\le c_t(x,m)
\le c_{t,0}^{n+1}(x,m)\le c_{t,0}^n(x,m) \le C \omega(m)$, where
$C$ is a
positive constant that depends solely on $\mu$. Note that
$p_u^{y,x}(m) \omega(m)$ is $m$-subadditive, which can be easily seen
by expanding $p_u^{y,x}$ as a series and using the subadditivity of
$a^{-1/2} \omega$ and the monotonicity of $a$. We next write
\begin{eqnarray*}
&&\biggl\langle\!\biggl\langle\frac{\omega(m')}{m'},
P_{t-s} [K_n(\tilde{\mu}_{s,0}^n) ] \biggr\rangle\!\biggr\rangle\\
&&\qquad=\frac{1}{2}
\int_{\mathbb{R}_+^2\times\mathbb{T}}
[ p_{t-s}^{y,x}(m'+m'') \omega(m'+m'')\\
&&\qquad\quad\hspace*{40.5pt}{}
-p_{t-s}^{y,x}(m') \omega(m')-p_{t-s}^{y,x}(m'') \omega(m'') ]\\
&&\hspace*{36.2pt}\qquad\quad{}
\times\frac{\kappa(m',m'')}{m'm''} 1_{E_n}(m'+m'') \tilde{\mu
}_{s,0}^n(y,dm')
\tilde{\mu}_{s,0}^n(y,dm'') \,dy.
\end{eqnarray*}
Together with
the subadditivity of $p_u^{y,x}(m') \omega(m')$,
the inequality $\tilde{\mu}_{s,0}^n\le\tilde{\mu}_{s,0}^{n+1} \le
\mu_s$ and the obvious
inclusion $E_n\subset E_{n+1}$, this implies
%
%
\begin{eqnarray}\label{equ57}
\omega(m) \biggl\langle\frac{\omega(m')}{m'},\nu^*\biggr\rangle
&\ge& \biggl\langle\!\biggl\langle\omega(m)\frac{\omega(m')}{m'},P_t \nu\biggl\rangle\!\biggl\rangle
\nonumber\\[-8pt]\\[-8pt]
&\ge&
c^n_{t,0}(x,m)\ge c^{n+1}_{t,0}(x,m)\ge c_t(x,m)
\ge0,\nonumber
\end{eqnarray}
as required.

We will need the fact that the solution $\mu_{\cdot}$ satisfies
%
%
\begin{eqnarray}\label{equ58}
\mu_t &=& \tilde{P}_t(c)\nu+\int_0^t\tilde{P}_{ts}(c)
[K^+(\mu _s) +
\delta_s \mu_s ] \,ds,\nonumber\\[-8pt]\\[-8pt]
\delta_t(x,m) &=& \biggl\langle\!\biggl\langle\frac{\omega(m) \omega(m')
-\kappa(m,m')}{m'},\mu_t\biggl\rangle\!\biggl\rangle.\nonumber
\end{eqnarray}
This can be proved by an argument similar to the one used to show the
equivalence of
(\ref{equ55}) and (\ref{equ51}).

Fix $k \in\mathbb{N}$, and assume that a sequence of measure paths and
mappings $\tilde{\mu}_{s,k}^n \le
\tilde{\mu}_{s,k}^{n+1}\le\mu_s$ and $c^n_{s,k}(x,m)
\ge c^{n+1}_{s,k}(x,m)\ge c_s(x,m)\ge0$ have already been defined. Let
\[
\delta^n_{t,k}=\biggl\langle\!\biggl\langle
\frac{\omega(m)\omega(m')-\kappa(m,m')}{m'},\tilde{\mu
}_{t,k}^n\biggl\rangle\!\biggl\rangle
\]
so that $0\le\delta^n_{t,k}\le\delta^n_{t,k+1}\le\delta_t$.
Finally, set
\begin{eqnarray*}
\tilde{\mu}^{n}_{t,k+1} &=& \tilde{P}_t(c^n_k)\nu^n +\int_0^t
\tilde{P}_{ts}(c^n_k) [K^+_n(\tilde{\mu}_{s,k}^n) + \delta
^n_{s,k}
\tilde{\mu}^n_{s,k} ] \,ds,\\
c^n_{t,k+1}(x,m) &=& \biggl\langle\!\biggl\langle\omega(m)\frac{\omega(m')}{m'},
P_t \nu\biggl\rangle\!\biggl\rangle\\
&&{} +\int_0^t \biggl\langle\!\biggl\langle\omega(m) \frac{\omega(m')}{m'},
P_{t-s} [K_n(\tilde{\mu}^n_{s,k+1}) ] \biggl\rangle\!\biggl\rangle \,ds.
\end{eqnarray*}
The Feynman--Kac formula (cf. \cite{RevuzYor}, Chapter 8) and identity
(\ref{equ58}) then imply
\[
0\le\tilde{\mu}^n_{s,k+1}\le\tilde{\mu}^{n+1}_{s,k+1} \le\mu_s
\]
and thus, by the same arguments that yield (\ref{equ57}),
\[
\omega(m) \biggl\langle\frac{\omega(m')}{m'},\nu^*\biggr\rangle\ge
c^n_{s,k+1}\ge c^{n+1}_{s,k+1}\ge c_s \ge0.
\]

The kernel $\tilde{\mu}^n_{s,k}$ is supported on $E_n$ for
each $n$ and all $k, s$. Also, $\tilde{\mu}^n_{s,k}\le\mu_s$
and the fact that $\omega(m)/m$ is bounded below in $E_n$
imply that the marginal density
$d\tilde{\mu}^n_{s,k}(dx,\mathbb{R}_+)/dx \in
L^{\infty}(dx)$ and its $L^{\infty}$-norm can be bounded uniformly in $k$,
%
%
\begin{equation}\label{equ59}
\biggl\| \frac{d\tilde{\mu}^n_{s,k}(dx,\mathbb{R}_+)}{dx} \biggr\|
_{\infty}
\le\beta_n \qquad\mbox{for all }
k\in\mathbb{N}.
\end{equation}

In view of these observations, we define a norm in
the vector space
\[
\biggl\{\eta\in\mathcaligr{M}_f(\mathbb{T}\times E_n), \frac
{d\eta
}{dx}(dx,\mathbb{R}_+)\in
L^{\infty}(dx) \biggr\}
\]
by
\[
\interleave\rho-\varrho\interleave= \bigl\| {\|\rho_x-\varrho_x\|}
\bigr\|_{\infty},
\]
if $\rho(dx,dm)=\rho_x(dm) \,dx$ and a similar formula holds for
$\varrho$.
Here $\|\rho_x-\varrho_x\|$ is the total variation norm of the signed
measure $\rho_x(dm)-\varrho_x(dm)$, which we may compute as
\[
\|\rho_x-\varrho_x\|=\int_{\mathbb{R}_+} |\rho_x-\varrho_x|(dm),
\]
if $|\rho_x-\varrho_x|(dm)$ denotes the total variation
of $\rho_x(dm)-\varrho_x(dm)$. Equivalently,
\[
\interleave\rho-\varrho\interleave
=\sup\biggl\{\langle f,\rho-\varrho\rangle| f\dvtx\int_{\mathbb
{T}}|\sup_m
f(x,m)| \,dx \le1 \biggr\}.
\]

Let $f(x,m)$ be such that $ \int_{\mathbb{T}}|\sup_m f(x,m)| \,dx
\le1 $.
We obtain
%
%
\begin{eqnarray}\label{equ510}
&&\langle f,\tilde{\mu}^n_{t,k+1}-\tilde{\mu}^n_{s,k} \rangle
\nonumber\\
&&\qquad=\langle
f,\tilde{P}_t(c_k^n)\nu^n-\tilde{P}_t(c_{k-1}^n)\nu^n \rangle
\nonumber\\[-8pt]\\[-8pt]
&&\qquad\quad{} + \int_0^t \langle f,
\tilde{P}_{ts}(c^n_k) [K_n^+(\tilde{\mu}_{s,k}^n) ]-
\tilde{P}_{ts}(c^n_{k-1}) [K_n^+(\tilde{\mu}_{s,k-1}^n) ]
\rangle \,ds
\nonumber\\
&&\qquad\quad{} + \int_0^t \langle f,
\tilde{P}_{ts}(c_k^n) [\delta^n_{s,k} \tilde{\mu}_{s,k}^n ]
-\tilde{P}_{ts}(c_{k-1}^n) [\delta^n_{s,k-1} \tilde{\mu
}_{s,k-1}^n] \rangle \,ds.\nonumber
\end{eqnarray}
We apply the Feynman--Kac formula to bound the first term on the right-hand side
of (\ref{equ510}) by
\[
\int_ {\mathbb{T}\times\mathbb{R}_+} E^{x,m} \biggl[ \biggl(\int_0^t
|c_{s,k}^n
(\chi_s,m)-c_{s,k-1}^n(\chi_s,m)| \,ds \biggr) f(\chi_t,m) \biggr] \nu
^n(dx,dm),
\]
where $\chi_s$ is a Brownian motion in $\mathbb{T}$ with diffusivity
$a(m)$. We
have
\begin{eqnarray*}
&&\int_0^t |c_{s,k}^n
(\chi_s,m)-c_{s,k-1}^n(\chi_s,m)| \,ds \\
&&\qquad\le\omega(m) \int_0^t \int_0^s
\biggl|\frac{\omega(m'+m'')}{m'+m''} p^{x',\chi_s}_{s-u}(m'+m'')
- \frac{\omega(m')}{m'} p^{x',\chi_s}_{s-u}(m') \biggr|\\
&&\qquad\quad\hspace*{54.6pt}{}\times\frac{\kappa(m',m'')}{m''}
1_{\{m'+m''\in E_n\}}\\
&&\qquad\quad\hspace*{54.6pt}{}\times|\tilde{\mu}_{u,k}^n(x',dm')\tilde{\mu}_{u,k}^n(x',dm'')
\\
&&\qquad\quad\hspace*{67.5pt}{}-
\tilde{\mu}_{u,k-1}^n(x',dm')\tilde{\mu}_{u,k-1}^n(x',dm'') |
\,dx' \,du \,ds.
\end{eqnarray*}
Note that in all these expressions the mass variable takes values
$m, m', m''$ or $m'+m''$ belonging to $E_n$, we may therefore replace
all functions that have it as an argument by an upper or lower bound,
as necessary, that do not depend on the mass variable.

Now, if $\rho$ and $\varrho$ are two finite, positive measures on
$(X,\Omega)$, then
$\rho\otimes\rho$ and $\varrho\otimes\varrho$ are finite,
positive measures on
$(X\times X,\Omega\times\Omega)$ and
\[
|\rho\otimes\rho- \varrho\otimes\varrho| \le
[\rho(X)+\varrho(X) ] |\rho-\varrho|.
\]
In particular,
\begin{eqnarray*}
&& |\tilde{\mu}_{u,k}^n(x',dm')\tilde{\mu}_{u,k}^n(x',dm'')-
\tilde{\mu}_{u,k-1}^n(x',dm')\tilde{\mu}_{u,k-1}^n(x',dm'') |\\
&&\qquad\le2 \beta_n
| \tilde{\mu}_{u,k}^n(x',dm')- \tilde{\mu
}_{u,k-1}^n(x',dm') |,
\end{eqnarray*}
$\beta_n$ as in (\ref{equ59}).

We thus have
\begin{eqnarray*}
&& \langle
f,\tilde{P}_t(c_k^n)\nu^n-\tilde{P}_t(c_{k-1}^n)\nu^n \rangle
\\
&&\qquad\le
\Gamma\int_0^t\int_0^s p(s-u,x',w) p(s,x,w) p(t-s,w,y) \\
&&\qquad\quad\hspace*{36.2pt}{} \times\biggl({\sup_m }|f(y,m)| \biggr)
\|\tilde{\mu}_{u,k}^n(x')-\tilde{\mu}_{u,k-1}^n(x') \|\\
&&\qquad\quad\hspace*{36.2pt}{} \times
\nu^n(dx,dm) \,dx' \,dw \,dy \,du \,ds,
\end{eqnarray*}
where $\Gamma$ is a positive constant that depends on $n$ and $\beta_n$
and can therefore be chosen uniformly in $k$.

Replace now
$\nu^n$ by its upper bound $\nu^*(dm) \,dx$, take the $L^{\infty}$ norm
of the total variation factor
$ |\tilde{\mu}_{u,k}^n(x')-\tilde{\mu}_{u,k-1}^n(x') |$,
and integrate in
$x, x', w$ and $s$, to obtain
\[
\langle
f,\tilde{P}_t(c_k^n)\nu^n-\tilde{P}_t(c_{k-1}^n)\nu^n \rangle
\le
\Gamma(T) \int_0^t
\interleave\tilde{\mu}_{u,k}^n-\tilde{\mu}_{u,k-1}^n\interleave
\,du.
\]
In this last step we used that $ \int_{\mathbb{T}}|{\sup_m f(x,m)}|
\,dx \le 1 $.

Similar computations yield
\begin{eqnarray*}
&&\int_0^t \langle f,
\tilde{P}_{ts}(c^n_k) [K_n^+(\tilde{\mu}_{s,k}^n) ]-
\tilde{P}_{ts}(c^n_{k-1}) [K_n^+(\tilde{\mu}_{s,k-1}^n) ]
\rangle
\,ds \\
&&\quad{} +\int_0^t \langle f,
\tilde{P}_{ts}(c_k^n) [\delta^n_{s,k} \tilde{\mu}_{s,k}^n ]
-\tilde{P}_{ts}(c_{k-1}^n) [\delta^n_{s,k-1} \tilde{\mu
}_{s,k-1}^n
] \rangle \,ds\\
&&\qquad\le\Gamma(T) \int_0^t
\interleave\tilde{\mu}_{u,k}^n-\tilde{\mu}_{u,k-1}^n\interleave
\,du
\end{eqnarray*}
and therefore, by (\ref{equ510}),
\[
\langle f,\tilde{\mu}^n_{t,k+1}-\tilde{\mu}^n_{s,k} \rangle\le
\Gamma(T) \int_0^t
\interleave\tilde{\mu}_{u,k}^n-\tilde{\mu}_{u,k-1}^n\interleave
\,du,
\]
$\Gamma(T)$ uniform in $k$, $f$. Take the
supremum over $f$ to conclude
that for $t\le T$,
%
%
\begin{equation}\label{equ511}
\interleave\tilde{\mu}_{t,k+1}^n-\tilde{\mu}_{t,k}^n\interleave
\le
\Gamma(T)
\int_0^t
\interleave\tilde{\mu}_{u,k}^n-\tilde{\mu}_{u,k-1}^n\interleave
\,du.
\end{equation}
Hence, by a standard contraction mapping argument, $\tilde{\mu
}^{t,n}_k$ converges
in $\mathcaligr{M}_f(\mathbb{T}\times E_n)$, uniformly in
$t\le T$.
The limit
is a continuous map $\tilde{\mu}^n\dvtx[0,T] \rightarrow\mathcaligr
{M}_f(\mathbb{T}
\times E_n)$
that satisfies (\ref{equ56}). Moreover, the properties that
$\tilde{\mu}^n_{s,k}\le\tilde{\mu}^{n+1}_{s,k} \le\mu_s$ and
$c_{s,k}^n\ge
c_{s,k}^{n+1} \ge c_s$ for
all $k, n \in\mathbb{N}^2 $ imply that the same holds for $\tilde
{\mu}^n$ and
$c^n$:
\[
\tilde{\mu}^n_s\le\tilde{\mu}^{n+1}_s \le\mu_s \quad\mbox
{and}\quad
c_s^{n}\ge c_s^{n+1}\ge c_s,
\]
for all $n\in\mathbb{N}$.

Suppose now that $\tilde{\mu}^n$ and $\tilde{\eta}^n$ are two
solutions in
$\mathcaligr{ D}(K)$ to
(\ref{equ56}) with respective values of the mapping $c$
defined as in the statement of
the lemma. A careful revision of the arguments leading to (\ref
{equ511}) shows
that the proof goes through verbatim if we replace $\tilde{\mu
}^n_{k+1}$ and
$\tilde{\mu}^n_k$ by fixed points of the iteration scheme $\tilde
{\mu}^n$ and
$\tilde{\eta}^n$. We thus have
\[
\interleave\tilde{\mu}_{t}^n-\tilde{\eta}^n_{t}\interleave\le
\Gamma(T)
\int_0^t
\interleave\tilde{\mu}_{u}^n-\tilde{\eta}_{u}^n\interleave \,du,
\qquad t\le T,
\]
which implies $\tilde{\mu}^n\equiv\tilde{\eta}^n$ by Gronwall's
lemma.
\end{pf}
\begin{pf*}{Proof of Theorem \protect\ref{theo3}}
Set $\lambda_0^n=1_{E_n^c}\nu$, and define the kernel $\tilde
{K}_n^-(\mu)=K^-(\mu) -K_n^-(\mu)$.
Given mappings $\tilde{\mu}^n_{\cdot}$ and $c^n_{\cdot}$ as in
Lemma \ref{lemma3},
define
\begin{eqnarray*}
\mu^n_t &=& P_t\nu^n + \int_0^t P_{t-s} [ K_n^+
(\tilde{\mu}_s^n)-K^-(\tilde{\mu}_s^n)-\gamma^n_s \tilde{\mu
}^n_s ] \,ds,\\
\lambda_t^n &=& P_t\lambda_0^n + \int_0^t
P_{t-s} [\tilde{K}_n^-(\tilde{\mu}^n_s)+\gamma_s^n
\tilde{\mu}_s^n ] \,ds,
\end{eqnarray*}
where
\[
\gamma_s^n(x,m)=c_s^n(x,m)-{\biggl\langle\!\biggl\langle\omega(m)\frac{\omega
(m')}{m'},\tilde{\mu}
_s^n \biggl\rangle\!\biggl\rangle}\ge c_s^n(x,m)-c_s(x,m)\ge0.
\]
The inequality in the last line follows from $\tilde{\mu}_s^n\le\mu
_s$ and
the definition of $c_s^n$; it has the consequence that $\lambda^n_s$
is a positive
measure for all $s\le T, n\in\mathbb{N}$.

We claim that $\mu^n=\tilde{\mu}^n$. Indeed, differentiating in the equations
satisfied by $\mu$ and $\tilde{\mu}$ shows that both
maps verify the equation
\[
\dot{\eta_t}=\frac{1}{2} a(m) \,\frac{\partial^2 }{\partial x^2}\,
\eta_t+K_n^+(\tilde{\mu}^n_t)-[c^n_t-\delta^n_t] \tilde{\mu}^n_t
\]
with initial value $\nu^n$, so that their difference is a weak solution
to $\dot{\eta}=\frac{1}{2} a(m) \,\frac{\partial^2 }{\partial x^2}\,
\eta_t$ started from the zero measure, and therefore it must be the
null measure. In particular this implies that $\mu_s^n\le\mu_s$.

We have, by the definition of $c^n_t$,
\begin{eqnarray*}
c_t^n(x,m) &=& {\biggl\langle\!\biggl\langle\omega(m)\frac{\omega(m')}{m'},P_t \nu\biggl\rangle\!\biggl\rangle}+\int
_0^t{\biggl\langle\!\biggl\langle\omega(m)\frac{\omega(m')}{m'},P_{t-s}
[K_n(\tilde{\mu}_s^n)]\biggl\rangle\!\biggl\rangle} \,ds\\
&=& {\biggl\langle\!\biggl\langle\omega(m)\frac{\omega(m')}{m'},\lambda_t^n\biggl\rangle\!\biggl\rangle}+{\biggl\langle\!\biggl\langle\omega
(m)\frac{\omega(m')}{m'},\mu_t^n\biggl\rangle\!\biggl\rangle}
\end{eqnarray*}
and on the other hand, from the definition of $\gamma_n^t$,
\[
c_t^n (x,m)=\gamma_t^n+\biggl\langle\!\biggl\langle\omega(m)\frac{\omega(m')}{m'},\tilde
{\mu}_t^n
\biggl\rangle\!\biggl\rangle.
\]
Hence
$\gamma_s^n=\llangle\omega(m)\frac{\omega(m')}{m'},\lambda_t^n\rrangle$.

From this point on, the proof is copied from that of
Theorem 5.4 in \cite{NorrisBcoag}. Let us define
$\alpha_t^n=\llangle\frac{\omega(m)}{m},\lambda_t^n\rrangle$. Then, due to the
subadditivity of $\omega(m)a(m)^{-1/2}$, we get
\[
\alpha_t^n + {\biggl\langle\!\biggl\langle\frac{\omega(m)}{m},
\mu_t^n \biggl\rangle\!\biggl\rangle} \ge \alpha_t^{n+1} + {\biggl\langle\!\biggl\langle\frac{\omega(m)}{m},
\mu_t^{n+1} \biggl\rangle\!\biggl\rangle}
\ge{\biggl\langle\!\biggl\langle\frac{\omega(m)}{m},\mu_t\biggl\rangle\!\biggl\rangle}
\]
for all $n\in\mathbb{N}$, and it follows from Lemma \ref{lemma3} that
$\alpha
^n_t\ge
\alpha_t^{n+1}$. We can hence define the monotone limits
\[
\alpha_t=\lim_{n\rightarrow\infty} \alpha_t^n,\qquad
\underline{\mu_t}=\lim_{n \rightarrow\infty} \mu_t^n,
\]
which satisfy
\[
\underline{\mu_t}\le\mu_t, \alpha_t+{\biggl\langle\!\biggl\langle\frac{\omega(m)}{m},
\underline{\mu}_t\biggl\rangle\!\biggl\rangle} \ge {\biggl\langle\!\biggl\langle\frac{\omega(m)}{m},
\mu_t \biggl\rangle\!\biggl\rangle}.
\]
Since $\omega(m)/m>0$, in the case that $\alpha_t=0$ a.e. we can
conclude
that $\underline{\mu_t}=\mu_t$ a.e., for all $t\in
[0,T]$. The uniqueness of the solution to (\ref{equ51}) will then follow
from the uniqueness of the solution to (\ref{equ56}).

We will now show that $\alpha_t$ vanishes. We start by proving that
$h_n(t)={\sup_{s\ge0}} \| \llangle\frac{\omega^2(m)}{m}, P_s\mu_t^n
\rrangle\|_{\infty}$ can be bounded uniformly in $n$ and $t\in[0,T]$.
We apply $P_s$ to the definition of
$\mu_t^n$, multiply by $\omega^2(m)/m$ and integrate over $E_n$ to
obtain
%
%
\begin{eqnarray}\label{equ512}
{\biggl\langle\!\biggl\langle\frac{\omega^2(m)}{m},P_s\mu_t^n \biggl\rangle\!\biggl\rangle}
&\le&
{\biggl\langle\!\biggl\langle\frac{\omega^2(m)}{m},
P_s \nu\biggl\rangle\!\biggl\rangle}\nonumber\\[-8pt]\\[-8pt]
&&{}+\int_0^t {\biggl\langle\!\biggl\langle\frac{\omega^2(m)}{m},
P_{s+t-r}K_n(\mu_r^n)\biggl\rangle\!\biggl\rangle} \,dr.\nonumber
\end{eqnarray}
By the subadditivity of $\omega$ and $\omega a^{-1/2}$, for any
$u\ge
0, x, z \in\mathbb{T}$, we have
\begin{eqnarray*}
&&\omega^2(m+m') p_u^{z,x}(m+m')-\omega^2(m) p_u^{z,x}(m)-\omega
^2(m') p_u^{z,x}(m')\\
&&\qquad\le
\omega(m) p_u^{z,x}(m) \omega(m')+\omega(m) \omega(m') p_u^{z,x}(m').
\end{eqnarray*}
Then
\begin{eqnarray*}
&&{\biggl\langle\!\biggl\langle\frac{\omega^2(m)}{m}, P_{s+t-r}K_n(\mu_r^n)\biggl\rangle\!\biggl\rangle}(x) \\
&&\qquad\le2 \int_{\mathbb{T}\times E_n^2}
\omega(m) p_{s+t-r}^{z,x}(m) \omega(m') \frac{\kappa
(m,m')}{mm'} \mu_{r}^n
(z,dm)\mu^n_r(z,dm') \,dz\\
&&\qquad\le2
\biggl\|{\biggl\langle\!\biggl\langle\frac{\omega^2}{m}
p_{s+t-r}^{\cdot,x},\mu^n_r\biggl\rangle\!\biggl\rangle}\,
{\biggl\langle\!\biggl\langle\frac{\omega \varpi}{m},\mu_r^n\biggl\rangle\!\biggl\rangle}\\
&&\qquad\quad\hspace*{9.1pt}{} + {\biggl\langle\!\biggl\langle\frac{\omega\varpi}{m}
p_{s+t-r}^{\cdot,x},\mu_r^n\biggl\rangle\!\biggl\rangle}\,{\biggl\langle\!\biggl\langle\frac{\omega^2}{m},\mu^n_r\biggl\rangle\!\biggl\rangle}
\biggr\|_1\\
&&\qquad\le2 \biggl\|{\biggl\langle\!\biggl\langle
\frac{\omega\varpi}{m},\mu^n_r\biggl\rangle\!\biggl\rangle}\biggr\|_{\infty}\,
{\biggl\langle\!\biggl\langle\frac{\omega^2}{m},P_{s+t-r}\mu^n_r\biggl\rangle\!\biggl\rangle}(x)\\
&&\qquad\quad{} + 2
\biggl\|{\biggl\langle\!\biggl\langle\frac{\omega^2}{m},\mu^n_r\biggl\rangle\!\biggl\rangle}\biggr\|_{\infty}\,
{\biggl\langle\!\biggl\langle
\frac{\omega\varpi}{m},P_{s+t-r}\mu^n_r\biggl\rangle\!\biggl\rangle}(x),
\end{eqnarray*}
where we used the bound $\kappa(m,m')\le\omega(m) \varpi(m')+
\varpi(m) \omega(m')$. Now, if we replace $\omega^2$ by
$\omega\varpi$ in (\ref{equ512}), the subadditivity of $a^{-1/2}
\omega\varpi$ implies that the time integral term is
nonpositive, and therefore
\[
\sup_r \sup_{s\ge
0} \biggl\|{\biggl\langle\!\biggl\langle\frac{\omega\varpi}{m},P_s\mu^n_r\biggl\rangle\!\biggl\rangle}\biggr\|_{\infty}
\le\biggl\langle\frac{\omega\varpi}{m},\nu^*\biggr\rangle<\infty
\]
by condition (\ref{equ54}) on $\nu$. This assumption also implies that
\[
\sup_{s\ge0} {\biggl\langle\!\biggl\langle\frac{\omega^2}{m},P_s\nu\biggl\rangle\!\biggl\rangle}\le\Gamma
<\infty
\]
for some positive constant $\Gamma$. Taking the supremum over $s\ge0$ in
(\ref{equ512}), we conclude that
\[
h_n(t)\le\Gamma+2\Gamma'\int_0^t h_n(r) \,dr\qquad \mbox{with
}\Gamma'>0.
\]
Note that the constants $\Gamma, \Gamma'$ can be chosen
independently of $n$. Then
$h_n(t)\le\Gamma e^{2\Gamma'T}$ holds uniformly in $n$ and $t\le T$, as
claimed.

We now consider the $L^1$ norm of $\alpha_t^n$. We replace
$\gamma_t^n(x,m)$ by its upper bound $\omega(m)
\alpha_t(x)$ in the definition of $\lambda^n$ and pass
to the limit as $n\rightarrow\infty$. By dominated convergence we have
\begin{eqnarray*}
\|\alpha_t\|_1 &\le& \int_0^t \int_{\mathbb{T}\times\mathbb
{T}\times\mathbb{R}_+}
\frac{\omega^2(m)}{m} \alpha_s(z) p_{t-s}^{z,x}(m)
\underline{\mu}_s(dm,dz) \,dx \,ds\\
&\le& \int_0^t \biggl(\sup_{n}\sup_{s\ge0} \biggl\|{\biggl\langle\!\biggl\langle\frac{\omega
^2}{m},\mu_s^n\biggl\rangle\!\biggl\rangle}\biggr\|_{\infty} \biggr)
\|\alpha_s\|_1 \,ds\\
&\le& \Gamma(T)\int_0^t \|\alpha_s\|_1 \,ds.
\end{eqnarray*}
Since $\|\alpha_t\|_1\le
\langle\omega/m,\nu^* \rangle<\infty$, this implies
$\alpha_t=0$ a.e., as required. The theorem follows.
\end{pf*}

\section*{Acknowledgments}
This work was done during the time I spent as a Research Associate
at the Statistical Laboratory, University of Cambridge.
I am grateful to James Norris for useful discussions and suggestions,
as well as for giving me the opportunity to read early versions of his
paper on Brownian coagulation \cite{NorrisBcoag}.

%

%
\printaddresses

\end{document}